\documentclass[openacc]{rsproca_new}
\usepackage{pifont}

\usepackage{algorithm}
\usepackage{algpseudocode}
\usepackage{amsmath}

\usepackage{xcolor}

\graphicspath{{RSPAInfo}{figure}}
\titlehead{Research}

\begin{document}

\title{WSINDy for Model Predictive Control with Applications to Fusion, Drones, and Chaos}

\author{
Cristian L\'{o}pez$^{1}$, 
Mckenna Partridge$^{1}$, Sebastian De Pascuale$^{2}$, Jeremy Lore$^{2}$, Andrew Christlieb$^{3}$, Stephen Becker$^{1}$, and David M.~Bortz$^{1}$}

\address{$^{1}$Department of Applied Mathematics, University of Colorado, Boulder, CO 80309-0526\\
$^{2}$Fusion Energy Division, Oak Ridge National Laboratory, 1 Bethel Valley Road,
Oak Ridge, TN 37830\\
$^{3}$Department of Computational Mathematics, Science and Engineering, Michigan State University, East Lansing, MI 48824}

\subject{applied mathematics, differential equations,  mathematical modeling}


\corres{Cristian L\'{o}pez, David M.~Bortz\\
\email{cristian.lopezruano@colorado.edu, david.bortz@colorado.edu}}

\begin{abstract}

\noindent The control of complex dynamical systems remains a fundamental challenge in science and engineering, where strong nonlinearities, the presence of noise and computational constraints often pose significant obstacles in traditional control approaches. Recent advances in data-driven methods, particularly system identification (SI) techniques, have shown a powerful alternative by providing fast, parsimonious, interpretable models that are well-suited for model predictive control (MPC). Building on these developments, the present article embeds weak sparse identification of nonlinear dynamics (WSINDy) with actuation control inputs (WSINDYc) within a MPC framework. Compared to benchmark data-driven methods, WSINDYc enables a more robust identification of the governing dynamics, particularly in the presence of high noise levels, resulting in more accurate and efficient control. The capabilities of the proposed WSINDY--MPC framework are demonstrated on a range of problems, including a tokamak plasma boundary model that includes main ion gas puff actuation, drone tracking and collision avoidance, the chaotic Lorenz system, and a simplified flight control model for an F-8 aircraft. The proposed framework achieves superior performance in the presence of noise, enabling longer prediction 
\absbreak 

\noindent horizons, lower trajectory tracking error, and a more reliable obstacle clearance, while simultaneously achieving lower MPC cost values compared to the baseline methods. Furthermore, we accelerate the optimization by compiling the MPC objective and constraints into C code, approaching real-time performance.

\end{abstract}

\rsbreak


\section{Introduction}

A central goal of modeling in the physical sciences is to identify mathematical relationships that accurately represent the behavior of real-world systems. The traditional approach for deriving governing equations from fundamental principles commonly uses conservation of momentum and energy. For systems whose behavior is well approximated as linear, identification can be performed using classical techniques such as the equation error model and linear regression \cite{Ljung1999}. However, modern systems often involve multiscale phenomena \cite{HaoYaoSuEtAl2024AdvNeuralInfProcessSyst}, nonlinear interactions \cite{Strogatz2015} and incomplete or unobservable information \cite{WiltingPriesemann2018NatCommun}. The mathematical complexity of these systems frequently exceeds our theoretical and computational ability to understand the model solutions. In recent years, the abundance of high resolution data has led to the development of data-driven methods \cite{BerkoozHolmesLumley1993AnnuRevFluidMech,LeeCarlberg2020JComputPhys}, rather than first principles; see \cite{Camps-VallsGerhardusNinadEtAl2023PhysicsReports} for a recent review of such methods. 

Over the past decades, numerous approaches have been developed for the identification of dynamical systems; notable works include an eigensystem realization algorithm \cite{JuangPappa1985JournalofGuidanceControlandDynamics}, symbolic regression \cite{SchmidtLipson2009Science}, Dynamic Mode Decomposition (DMD) \cite{Schmid2010JFluidMech,TuRowleyLuchtenburgEtAl2014JComputDyn}, neural networks (NN)\cite{ChenRubanovaBettencourtEtAl2018AdvancesinNeuralInformationProcessingSystems,RaissiKarniadakis2018JournalofComputationalPhysics,LiuMaWangEtAl2024arXiv240810205}, spectral submanifolds \cite{CenedeseAxasBauerleinEtAl2022NatCommun}, and regression on prescribed function bases \cite{CrutchfieldMcNamara1987ComplexSyst,WangYangLaiEtAl2011PhysRevLett}. From a model-selection perspective, Akaike's information-theoretic criterion \cite{Akaike1974IEEETransAutomControl} is a commonly used metric to determine the trade-off between model fit and complexity. Building on these ideas, in \cite{BruntonProctorKutz2016ProcNatlAcadSci}, the sparse identification of nonlinear dynamics (SINDy) method was proposed as a sparsity-promoting model selection framework for discovering governing equations from data. Subsequently, this approach was extended in \cite{MessengerBortz2021MultiscaleModelSimul} by using a weak formulation, preserving the model selection structure while enhancing robustness to noise.

Among these techniques, the SINDy method \cite{BruntonProctorKutz2016ProcNatlAcadSci} has gained considerable attention for discovering the governing equations of dynamical systems from data. It assumes the system's dynamics can be represented as a linear combination of a few terms selected from a large library of candidate functions, such as polynomials, trigonometric, rational terms, etc., that might appear in the governing equations. After collecting sample measurements and estimating their time derivatives, often through finite differences or total variation regularization to manage noise, the algorithm then performs a sequentially thresholded, least-squares regression between the numerical derivatives of the data and the evaluated library of functions to identify/estimate the coefficients of the terms that represent the data. Due to its computational efficiency, interpretable results, and accuracy under low noise levels, SINDy has found applications across a wide range of domains \cite{ManganBruntonProctorEtAl2016IEEETransMolBiolMulti-ScaleCommun,SorokinaSygletosTuritsyn2016OptExpress,BoninsegnaNuskeClementi2018JChemPhys,LaiNagarajaiah2019MechSystSignalProcess,JiangXiongZhangEtAl2021NonlinearDyn,FukamiMurataZhangEtAl2021JFluidMech}. Beyond system identification of uncontrolled systems, SINDy has also been extended to problems in control (SINDYc) \cite{BruntonProctorKutz2016IFAC-PapersOnLine},  where the identified dynamics accounts for the influence of actuation inputs. Once such a model is obtained, it can be embedded into advanced control strategies. 

In control theory, several strategies have been developed to effectively regulate dynamical systems. For instance, proportional-integral-derivative controllers \cite{BoraseMaghadeSondkarEtAl2021IntJDynamControl} are widely accepted owing to their simplicity and ease of implementation; linear quadratic regulators \cite{BruntonKutz2019} achieve optimal control with an analytical solution under a quadratic cost function for linear systems; NNs have demonstrated the ability to model and control complex systems~\cite{Bottcher2026NonlinearDyn}, including recurrent and ensemble architectures , when sufficient training data and computational resources are available; and model-free adaptive control methods construct control laws directly from online input-output data\cite{HeydariNovinzadehTayefi2025AdaptiveControlampSignal}. Another prominent strategy that has gained significant attention is MPC, which is the focus of our work. This model-based optimization framework computes control actions over a receding horizon and continuously updates its course using system feedback from measurements \cite{GarciaPrettMorari1989Automatica, SchwenzerAyBergsEtAl2021IntJAdvManufTechnol}. It has become a cornerstone for controlling nonlinear systems with instabilities, constraints, and time delays \cite{BruntonZolmanKutzEtAl2025AnnuRevControlRobotAutonSyst}. MPC strongly depends on model quality to calculate optimal control actions, as it uses this model to predict the system's future behavior over the horizon and then selects the best current action based on those predictions. A broad range of techniques for obtaining data-driven models has therefore been integrated with MPC, including Gaussian process\cite{CaoLaiAlam2017JIntellRobotSyst}, DMD with control \cite{ProctorBruntonKutz2016SIAMJApplDynSyst}, SINDYc \cite{KaiserKutzBrunton2018ProcRSocA,FaselKaiserKutzEtAl2021202160thIEEEConfDecisControlCDC}, Koopman operator \cite{KordaMezic2018Automatica}, recurrent NNs \cite{WuTranRinconEtAl2019AIChEJournal,WuTranRinconEtAl2019AIChEJournala,WuChristofidesWuEtAl2025RevChemEng}.

In particular, the interpretable models provided by SINDYc have been combined with MPC (SINDY--MPC) \cite{KaiserKutzBrunton2018ProcRSocA,FaselKaiserKutzEtAl2021202160thIEEEConfDecisControlCDC} to control strongly nonlinear systems, demonstrating superior performance in both prediction and control compared to techniques such as DMDc and NNs. 
Additional applications of the SINDY--MPC methodology include chemical processes \cite{AbdullahWuChristofides2021ComputChemEng,AbdullahAlhajeriChristofides2022IndEngChemRes, AbdullahChristofides2023ComputChemEng,AbdullahChristofides2023ChemicalEngineeringResearchandDesign}, wearable pneumatic systems \cite{LeeRenQianEtAl2025IEEEASMETransMechatron}, diesel engine air path system \cite{YahagiSetoYonezawaEtAl2025IntJControlAutomSyst}, unmanned aerial vehicles \cite{GuevaraVarela-AldasGandolfoEtAl2025Drones,LeeKimKimEtAl2025IETControlTheoryAppl}, and plasma boundary control \cite{LoreDePascualeLaiuEtAl2023NuclFusion}. Although successful in the tested systems, the traditional SINDy method relies on pointwise derivatives of the data, which reduces its robustness to noise and potentially hinders its applicability in MPC scenarios. To address this limitation, several approaches have been proposed in recent years that replace explicit derivatives with integrals in the sparse model selection \cite{SchaefferMcCalla2017PhysRevE,GurevichReinboldGrigoriev2019Chaos,ReinboldGrigoriev2019PhysRevE,MessengerBortz2021MultiscaleModelSimul,MessengerBortz2021JComputPhys,MessengerBortz2022PhysicaD,BortzMessengerDukic2023BullMathBiol,BortzMessengerTran2024NumericalAnalysisMeetsMachineLearning,TranBortz2025arXiv250703206}, leading to a weak form of SINDy. As shown in \cite{MessengerBortz2021MultiscaleModelSimul},\cite{MessengerBortz2021JComputPhys}, and \cite{LopezNaranjoSalazarEtAl2025JComputPhys}, this formulation is computationally efficient due to its convolutional-based implementation and also improves robustness to noise, since derivatives act on smooth, localized test functions, providing a natural filtering effect that reduces the effect of noise. Such weak-form constructions date back to 1954 \cite{Shinbrot1954NACATN3288} for parameter estimation in differential equations.

Weak SINDy (WSINDy) has been utilized in numerous fields, including population dynamics \cite{MessengerDwyerDukic2024JRSocInterface}, weather \cite{MinorMessengerDukicEtAl2025JournalofGeophysicalResearchMachineLearningandComputation}, magnetohydrodynamics \cite{VaseyMessengerBortzEtAl2025JComputPhys}, mechanical structures \cite{SchmidDoostanPourahmadian2024arXiv240920510}, and resonance imaging \cite{WoodallEsparzaGutovaEtAl2024APLBioeng}, to name a few. In this work, we start with the method of WSINDy with actuation control inputs (WSINDYc) to identify interpretable nonlinear models that incorporate the effects of actuation inputs \cite{LopezNaranjoSalazarEtAl2025JComputPhys}. We then integrate it with MPC to enable data-driven closed-loop control (WSINDY--MPC) of a range of nonlinear systems. 
We show that WSINDY--MPC outperforms relevant baseline methods (SINDYc, E-SINDYc, DMDc and a NN approach) in terms of noise robustness, prediction and control performance. We also introduce an ensemble version, E-WSINDY--MPC, using the ensembling technique of \cite{FaselKutzBruntonEtAl2022ProcRSocMathPhysEngSci}, and see that WSINDY--MPC and E-WSINDY--MPC have comparable performance,  with neither method consistently dominating across the tested cases and evaluation criteria, though the ensemble method is computationally slower.

The remainder of this paper is organized as follows. Sections \ref{sec:Sparse Modeling} and \ref{sec:MPC} present the theory of sparse nonlinear modeling and MPC, respectively. Section \ref{sec:WSINDY-MPC} introduces the proposed WSINDY--MPC framework. Section \ref{sec:Comparison} provides the benchmark comparison across a variety of systems, namely fusion, drone, and Lorenz attractor. Finally, section \ref{sec:Conc} contains the concluding remarks and outlines potential directions for future research.   

\section{Sparse Nonlinear Modeling}
\label{sec:Sparse Modeling}

Consider a dynamical system in $D$ dimensions governed by first-order equations of the form
\begin{align}\label{2.1}
\dot{\mathbf{x}}(t) = \mathbf{f}(\mathbf{x}(t)), \quad 
\mathbf{x}(0) = \mathbf{x}_0 \in \mathbb{R}^D, \quad 0 \le t \le T_t,
\end{align}
where $\mathbf{x}(t) = [x_1(t), x_2(t), \dots, x_D(t)]^\top \in \mathbb{R}^D$ denotes the system state, $\mathbf{f}$ is in general a nonlinear function that describes the system dynamics, and $T_t$ is the total time.  

The SINDy algorithm \cite{BruntonProctorKutz2016ProcNatlAcadSci} approximates the derivative $\dot{\mathbf{x}}(t)$ as a linear combination of a library of candidate functions $f_j$, $j=1,\dots,J$, that may describe the system dynamics. These candidate functions are collected into a library $\boldsymbol{\Theta}(\mathbf{x}) = [f_1(\mathbf{x}),\ f_2(\mathbf{x}),\ \dots,\ f_J(\mathbf{x})] \in \mathbb{R}^{1\times J}$, where each $f_j(\mathbf{x})$ denotes a pointwise evaluation.  Given measurements at discrete time instants $t_k$, $k=1,\dots,N$, the library evaluated on the data $
\mathbf{X} = [\mathbf{x}(t_1);\ \mathbf{x}(t_2);\ \dots;\ \mathbf{x}(t_N)] \in \mathbb{R}^{N\times D}$ is denoted by $\boldsymbol{\Theta}(\mathbf{X}) \in \mathbb{R}^{N\times J}$.

To identify the fewest terms that best capture the dynamics, a sparse regression problem is formulated as
\begin{align}\label{2.2}
\mathbf{W} = \arg \min_{\mathbf{W}} 
\|\dot{\mathbf{X}} - \boldsymbol{\Theta}(\mathbf{X})\mathbf{W}\|_2^2 
+ \lambda \|\mathbf{W}\|_0,
\end{align}
where $\dot{\mathbf{X}} \in \mathbb{R}^{N\times D}$ are the estimated time derivatives of the data,  $\mathbf{W} = [\mathbf{w}_1,\ \mathbf{w}_2,\ \dots,\ \mathbf{w}_D] \in \mathbb{R}^{J\times D}$ is the coefficient matrix, with each $\mathbf{w}_d \in \mathbb{R}^J$ corresponding to the dynamics of the $d$th state variable, $\|\mathbf{W}\|_0$ is the $\ell_0$ pseudo-norm which is the count of the number of nonzero entries of the matrix, and the hyperparameter $\lambda>0$ controls the sparsity of the solution. While solving Eq.~\eqref{2.2} exactly is generally intractable, 
there are efficient algorithms to approximately solve it. In particular, the problem is separable in the columns of $\mathbf{W}$, and 
each column $\mathbf{w}_d$ is typically found using a sequential thresholding least squares method.

After solving Eq.~\eqref{2.2}, the model $\mathbf{f}$ is approximated as
\begin{align}\label{2.3}
\dot{\mathbf{x}}(t) = \mathbf{f}(\mathbf{x}(t)) 
\approx \boldsymbol{\Theta}(\mathbf{x}(t))\mathbf{W}.
\end{align}

\subsection{Sparse identification of nonlinear dynamics with control}\label{sec:SINDYc}
The SINDYc framework \cite{BruntonProctorKutz2016IFAC-PapersOnLine} was introduced shortly after the original SINDy algorithm, generalizing the method to account for the effect of actuation inputs. In this case, the dynamics are modeled as 
\begin{align}\label{2.4}
\dot{\mathbf{x}}(t) = \mathbf{f}(\mathbf{x}(t), \mathbf{u}(t)).
\end{align}
with the corresponding discrete-time representation
\begin{align}\label{2.5}
\mathbf{x}_{k+1} = \mathbf{F}(\mathbf{x}_k, \mathbf{u}_k),
\end{align}
where $\mathbf{F}$ is the exact state update map.

Accordingly, the library of candidate functions is extended to include functions of both the state and the input, $\boldsymbol{\Theta}(\mathbf{x},\mathbf{u}) = [f_1(\mathbf{x},\mathbf{u}),\ f_2(\mathbf{x},\mathbf{u}),\ \dots \ , \ f_J(\mathbf{x},\mathbf{u})] \in \mathbb{R}^{1\times J}$, where each $f_j(\mathbf{x},\mathbf{u})$ denotes a candidate term that may describe the controlled dynamics. Typical choices include $\boldsymbol{\Theta(\mathbf{X,U})}=[\mathbf{1},\ \mathbf{X},\ \mathbf{U},\ (\mathbf{X}\otimes\mathbf{X}),\ (\mathbf{X}\otimes\mathbf{U}),\ (\mathbf{U}\otimes\mathbf{U}),\ \dots \ ]$, where $\mathbf{X}\otimes\mathbf{U}$ denotes the matrix consisting of all product combinations of the components in $\mathbf{x}$ and $\mathbf{u}$, and $\mathbf{U}=[\mathbf{u}(t_1);\mathbf{u}(t_2);\ \dots \ ; \ \mathbf{u}(t_N)]\in\mathbb{R}^{N\times V}$, with $V$ denoting the number of control inputs. With this construction,  Eq.~\eqref{2.4} can be written in the regression form 
\begin{align}\label{2.6}
\dot{\mathbf{x}}(t) = \boldsymbol{\Theta}(\mathbf{x}(t),\mathbf{u}(t))\mathbf{W}.
\end{align}

The estimation of $\mathbf{W}$ follows the same methodology as that used in Eq.~\eqref{2.2}.

\subsection{Weak sparse identification of nonlinear dynamics with control}\label{sec:WSINDYc}
When working with noisy data, estimating derivatives is difficult, as small fluctuations in the data can create large errors. This poses a major challenge in the SINDy method \cite{BruntonProctorKutz2016ProcNatlAcadSci} and its extensions \cite{BruntonProctorKutz2016IFAC-PapersOnLine,KaiserKutzBrunton2018ProcRSocA,FaselKutzBruntonEtAl2022ProcRSocMathPhysEngSci}, which depend on their estimation. The WSINDy method addresses this by rewriting Eq.~\eqref{2.1} in its integral form. Instead of computing derivatives, WSINDy integrates the data against selected test functions and, through integration by parts, transfers derivatives from the data to these test functions. This approach smooths out noise and naturally filters out high-frequency components, leading to more robust system identification without the instability caused by numerical differentiation.

Recently, we also extended the WSINDy framework for ordinary differential equations (ODEs) to incorporate actuation inputs \cite{MessengerBortz2021JComputPhys, LopezNaranjoSalazarEtAl2025JComputPhys}, allowing the weak-form approach to handle externally actuated systems. 
Related weak-form sparse identification ideas have been used in control frameworks, including micro autonomous surface vehicles \cite{ChenWang2025arXiv250906882} and virtual coupling trains \cite{WangLiWenEtAl2026IEEETransVehTechnol}. The former was developed on the earliest version of WSINDy \cite{MessengerBortz2021MultiscaleModelSimul}, where the test function's parameters depend on manually defined values, while the latter uses the integral form with a structured pruning. By contrast, the present work focuses on using test functions whose parameters are automatically determined from data \cite{MessengerBortz2021JComputPhys, LopezNaranjoSalazarEtAl2025JComputPhys}.

The WSINDYc framework \cite{LopezNaranjoSalazarEtAl2025JComputPhys} is obtained by converting Eq.~\eqref{2.4} to its weak form. This is done by multiplying Eq.~\eqref{2.5} by a compactly supported test function $\phi(t)\in C^1_c([0,T_t])$, with $\mathrm{supp}(\phi)\subset(0,T_t)$, and integrating by parts to obtain
\begin{align}\label{2.7}
-\mathit{\int}_0^{T_t} \dot{\phi}(t) \, \mathbf{x}(t) \, dt
= \mathit{\int}_0^{T_t} \phi \, \boldsymbol{\Theta}(\mathbf{x}(t), \mathbf{u}(t)) \, \mathbf{W} \, dt.
\end{align}

By using translations of a fixed $\phi(t)$ through the time domain, the weak form integrals in Eq.~\eqref{2.6} take the form of discrete convolutions, which are then approximated using the trapezoidal quadrature scheme
\begin{align}\label{2.8}
-\boldsymbol{\dot{\phi}}(t)\mathbf{X}(t) \approx \boldsymbol{\phi}(t)\boldsymbol{\Theta}(\mathbf{X}(t),\mathbf{U}(t))\mathbf{W},
\end{align}
where $\boldsymbol{\phi}(t)$, $\boldsymbol{\dot{\phi}}(t)\in\mathbb{R}^{1\times N}$ denote $\boldsymbol{\phi}(t):=[ \phi(t_1)\ \ ...\ \phi(t_N) ] \mathcal{Q}$ and $\boldsymbol{\dot{\phi}}(t):=[ \dot{\phi}(t_1)\ \ ...\ \dot{\phi}(t_N) ] \mathcal{Q}$, respectively, with $\mathcal{Q}\in\mathbb{R}^{1\times N}$ the quadrature vector. Thus, from Eq.~\eqref{2.7} we obtain the weak convolution equation 
\begin{align}\label{2.9}
\mathbf{B} \approx \mathbf{GW},
\end{align}
\noindent where $\mathbf{B}\in\mathbb{R}^{2N-1}$ and $\mathbf{G}\in\mathbb{R}^{(2N-1) \times J}$. Note that both $\boldsymbol{\dot{\phi}}(t)$ and $\mathbf{X}(t)$ are sampled at $N$ time points; their full discrete convolution has length $2N-1$. However, since $\boldsymbol{\dot{\phi}}(t)$ is localized, with length $2m+1$ time points ($m$ being its length support), in practice, we retain only the interior values, corresponding to the complete overlap between the data and the test function, which removes boundary effects and yields $N-2m$ entries. Additionally, for any dynamical system, we use the corresponding $\boldsymbol{\phi}(t)$ and $\mathbf{G}$ matrix for each component. Similar to Eq.~\eqref{2.2}, we seek a sparse solution to the least-squares problem
\begin{align}\label{2.10}
\mathbf{W} \approx \arg \min_{\mathbf{W}} 
\|\mathbf{B} - \mathbf{GW}\|_2^2 
+ \lambda \|\mathbf{W}\|_0.
\end{align}
We approximate its solution using the modified sequential thresholding least squares (MSTLS) algorithm \cite{MessengerBortz2022PhysicaD,LopezNaranjoSalazarEtAl2025JComputPhys}. MSTLS is a greedy method that alternates between least-squares regression and coefficient thresholding. The threshold parameter $\lambda$ is selected from a grid of $100$ candidate values, $\{10^{-4+\frac{4\ell}{99}} \mid \ell\in\{0,...,99\}\}$, by  choosing the one that minimizes the loss function
\begin{align}\label{2.11}
\mathcal{L}=\frac{\|\mathbf{G}(\mathbf{W^{\lambda}-W^{0}})\|_2}{\|\mathbf{GW}^{0}\|_2} + \frac{\#\{\mathcal{I}^{\lambda}\}}{J},
\end{align}

\noindent where $\mathbf{W}^{0}$ is the least squares solution of $\mathbf{GW}=\mathbf{B}$, $\#\{\cdot\}$ represents cardinality, and $\{\mathcal{I}^{\lambda}\}$ is the index set of non-zeros coefficients of $\mathbf{W}^{\lambda}$.

\section{Model Predictive Control}
\label{sec:MPC}

MPC is a well-established, model-based control technique that has been adopted across various industrial applications due to its proven effectiveness \cite{GarciaPrettMorari1989Automatica,SchwenzerAyBergsEtAl2021IntJAdvManufTechnol,BiekerPeitzBruntonEtAl2020TheorComputFluidDyn,BruntonZolmanKutzEtAl2025AnnuRevControlRobotAutonSyst}. If a precise model and adequate computing resources are available, MPC provides a versatile framework for controlling strongly nonlinear systems with constraints, uncertainties, and time delays \cite{KaiserKutzBrunton2018ProcRSocA,FaselKutzBruntonEtAl2022ProcRSocMathPhysEngSci}.

In MPC, given the system state $\mathbf{x}_j$, a control sequence $\mathbf{u}(\mathbf{x}_j) = \{\mathbf{u}_{j+1}, \ \dots, \  \mathbf{u}_{j+m_c}\}$ is obtained by solving an optimization problem over a receding horizon $T_c = m_c\Delta t^\text{M}$, with $m_c$ as the number of control steps and $\Delta t^\text{M}$ is the internal model integration step, which may differ from the measurements sampling time. To obtain the optimal control sequence $\mathbf{u}(\mathbf{x}_j)$, the cost function
\begin{align}\label{3.1}
J = \sum_{k=0}^{m_p-1} \left\| \hat{\mathbf{x}}_{j+k} - \mathbf{r}_k \right\|_{\mathbf{Q}}^{2} 
    + \sum_{k=1}^{m_c-1} \Big(\left\| \mathbf{u}_{j+k} \right\|_{\mathbf{R_u}}^{2} + \left\| \mathbf{\Delta u}_{j+k} \right\|_{\mathbf{R}_{{\Delta} u}}^{2}\Big),
\end{align}
is minimized over a prediction horizon $T_p = m_p\Delta t$, with $m_p$ as the number of prediction steps. Generally $T_c \le T_p$; if $T_c < T_p$, the input $\mathbf{u}$ is kept fixed after $T_c$, i.e., $\mathbf{u}_{j+k}=\mathbf{u}_{j+m_c}$, $k>m_c$. Here, $ \hat{\mathbf{x}}$ denotes the predicted system state, initialized with the measured state $\hat{\mathbf{x}}_j=\mathbf{x}_j$, $\mathbf{r}$ the reference state to be tracked, $\left\| \mathbf{x} \right\|_{\mathbf{Q}}^{2} = \mathbf{x}^\top \mathbf{Q}\mathbf{x}$, $\Delta \mathbf{u}_k = \mathbf{u}_k - \mathbf{u}_{k-1}$, and $\mathbf{Q} \succeq 0$ (positive semi-definite), and $\mathbf{R}_u \succ 0$, $\mathbf{R}_{\Delta u} \succ 0$ (positive definite) are weight matrices. These matrices encode the trade-off between tracking accuracy and control effort: larger entries $\mathbf{Q}$ prioritize state regulation, while larger entries in $\mathbf{R}_u$ and $\mathbf{R}_{\Delta u}$ penalize actuation magnitude and rate of change, respectively, yielding smoother but slower responses. In practice, these matrices are selected through iterative tuning \cite{RawlingsMayneDiehl2024}. The optimization problem is subject to the discrete-time prediction model 
\begin{align}\label{3.2}
\hat{\mathbf{x}}_{k+1} =  
    \hat{\mathbf{F}}(\hat{\mathbf{x}}_{k},\mathbf{u}_{k}),
\end{align}
where $\hat{\mathbf{F}}$ denotes the state update map obtained by numerically integrating Eq.~\eqref{2.4} over one sampling interval (i.e., using a forward 
fourth-order Runge-Kutta (RK4) scheme), and the constraints:
\begin{align}
\mathbf{u}_{\text{min}} &\le \mathbf{u}_k \le \mathbf{u}_{\text{max}}, \label{3.3}\\
\mathbf{\Delta u}_{\text{min}} &\le \mathbf{\Delta u}_k \le \mathbf{\Delta u}_{\text{max}}. \label{3.4}
\end{align}

At each \emph{measurement sampling} time step $\mathbf{u}(\mathbf{x}_j)$ is calculated, however, only the first prescribed control input $\mathbf{u}_{j+1}$ is physically applied to the system. 
The true plant dynamics $\dot{\mathbf{x}} = \mathbf{f}_{\mathrm{plant}}(\mathbf{x},\mathbf{u})$ are then propagated over the interval $T_s$ to obtain the next state $\mathbf{x}_{j+1}$. In our simulation setting, this evolution is approximated via RK4 with time step $\Delta t^{\text{plant}}$ and integration substeps $n_s^\text{plant}$. 

Unlike \cite{KaiserKutzBrunton2018ProcRSocA}, where noiseless state feedback is assumed, the state fed back to the controller is corrupted by noise, $\mathbf{y}_{j+1} = \mathbf{x}_{j+1} + \boldsymbol{\varepsilon}_{j+1}$, with $\boldsymbol{\varepsilon}\in \mathcal{N}(0, \sigma^{2})$. The noise magnitude is defined as $\eta_i = \sigma / \text{std}(x_i)$, and the states provided to the controller are corrupted using the same noise level as that present in the system identification data. The $\text{std}(x_i)$ represents the sample standard deviation of the $i$-th state variable over the training data. Accordingly, the noisy training data is denoted $\mathbf{Y} = \mathbf{X}+\mathbf{E} \in \mathbb{R}^{N\times D}$, where  $\mathbf{E}$ is the corresponding noise matrix. 

The optimization is reinitialized with the updated state $\mathbf{y}_{j+1}$, and a new sequence $\mathbf{u}(\mathbf{y}_{j+1})$ is computed. As in \cite{KaiserKutzBrunton2018ProcRSocA}, for a fair comparison, all methods use the same optimization approach through MATLAB's \texttt{fmincon} function, configured to employ the sequential-quadratic-programming algorithm and limited to $100$ iterations.

\section{WSINDy for model predictive control}
\label{sec:WSINDY-MPC}

\begin{figure}[tb]
\centering\includegraphics[width=5.0in]{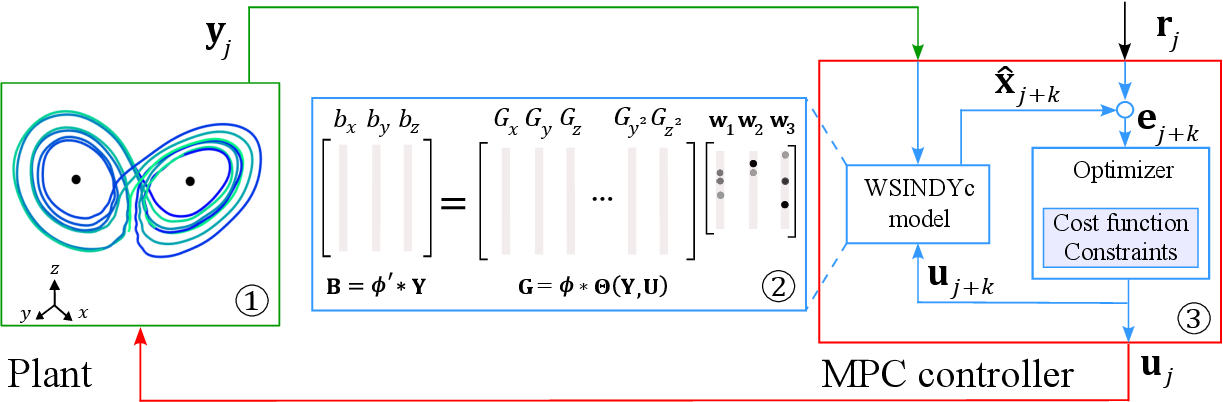}
\caption{The WSINDY--MPC framework. Offline stage: \ding{192}~data generation and \ding{193}~model identification using WSINDYc. Online stage: \ding{194}~MPC implementation using the learned model. For example, the Lorenz 63 system, shown in plant, is driven toward a selected fixed point (marked as black dots) using MPC.}
\label{MPC_flowchart}
\end{figure}

Since the MPC depends on the accuracy of a mathematical model, the use of noise-robust identification methods is essential for reliable control. Usually, SINDy-based approaches rely on numerical differentiation of time-series data, which is highly sensitive to noise. Various strategies have therefore been proposed to improve robustness, including smoothing, filtering techniques, regularized differentiation methods \cite{BruntonProctorKutz2016ProcNatlAcadSci,LoreDePascualeLaiuEtAl2023NuclFusion,LejarzaBaldea2022SciRep,WentzDoostan2023ComputerMethodsinAppliedMechanicsandEngineering}, etc. However, these approaches often require careful parameter tuning, may suppress important dynamical features, and might increase computational cost. By contrast, WSINDYc \cite{MessengerBortz2021MultiscaleModelSimul,MessengerBortz2021JComputPhys,LopezNaranjoSalazarEtAl2025JComputPhys} adopts a weak formulation in which numerical differentiation of the data is avoided by integrating against smooth test functions, by transferring the derivatives to the test functions. From a signal processing standpoint, this operation acts as a low-pass filter on the data, attenuating high-frequency noise and therefore improving the system identification robustness to noisy measurements. This attribute motivates its integration with MPC for data-driven control of strongly nonlinear systems, particularly when the collected measurement data are corrupted by significant noise. The proposed WSINDY--MPC framework is depicted in Fig.~\ref{MPC_flowchart} and consists of two stages: an offline stage, in which informative data are collected and an ODE model is identified from these data using WSINDYc, and an online stage, where MPC is implemented using the learned model to compute control actions. We illustrate the pseudocode of the WSINDY--MPC framework in  Algorithm ~\ref{alg:WSINDY_MPC}.

\begin{table}[t]
\centering
\refstepcounter{table}
\label{alg:WSINDY_MPC}

\parbox{\columnwidth}{%
\hrule
\vspace{3pt}
\textbf{Algorithm 1. WSINDY--MPC framework}
\vspace{3pt}
\hrule
}

{\ttfamily
\parbox{\columnwidth}{
\begin{tabular}{@{}p{0.1\columnwidth} p{0.85\columnwidth}}
\textbf{Input:} &
\emph{System Identification (SI)}: Training data $\mathbf{Y},\mathbf{U}$ (sampled at $\Delta t$),
candidate library $\boldsymbol{\Theta}(\cdot)$, weak-form settings (test function type, query points); \\
&
\emph{MPC}: horizon steps $m_p$ and $m_c$,update interval $T_s$,
model time step $\Delta t^{\mathrm{M}}$, model substeps $n_s^{\mathrm{M}}$, weights $\mathbf{Q},\mathbf{R}_u,
\mathbf{R}_{\Delta u}$, constraints $\Delta \mathbf{u}_{\min},
\Delta \mathbf{u}_{\max}, \mathbf{u}_{\min}, \mathbf{u}_{\max}$,
reference $\mathbf{r}$. \\[6pt]
\textbf{Output:} &
Sequence of control inputs $\mathbf{u}_j$. \\
\end{tabular}
}
}

\vspace{6pt}

{\ttfamily
\parbox{\columnwidth}{
\begin{tabular}{@{} r @{\quad} p{0.92\columnwidth} @{}}
\multicolumn{2}{@{}l@{}}{\emph{Offline}: SI via WSINDYc}\\[6pt]

1: & Create $\mathbf{B},\mathbf{G}$ 
   \hfill $\triangleright$ Weak regression matrices, see \cite{MessengerBortz2021JComputPhys,LopezNaranjoSalazarEtAl2025JComputPhys} \\

2: & $\displaystyle\mathbf{W} \approx \operatorname*{\arg\min}_{\mathbf{W}}
   \|{\mathbf{B}} - \mathbf{GW}\|_2^2 + \lambda \|\mathbf{W}\|_0$ 
   \hfill $\triangleright$ Solve Eq.~\eqref{2.9} using $\mathrm{MSTLS}$ \\

3: & $\hat{\mathbf{f}}(\mathbf{y},\mathbf{u})
   = \boldsymbol{\Theta}(\mathbf{y},\mathbf{u})\,\mathbf{W}$
   \hfill $\triangleright$ Obtain learned continuous-time model \\

4: & $n_s^{\mathrm{M}} = \operatorname{round}\!\left({T_s}/{\Delta t^{\mathrm{M}}}\right)$
   \hfill $\triangleright$ Get substeps \\

5: & $\hat{\mathbf{F}} \leftarrow \texttt{RK4} \!\left(\hat{\mathbf{f}},\Delta t^{\mathrm{M}},n_s^{\mathrm{M}}\right)$
   \hfill 
   \begin{tabular}[t]{@{}r@{}}
   $\triangleright$ Define forward operator \\
   (single-step or substepped RK4)
   \end{tabular} \\[8pt]

\multicolumn{2}{@{}l@{}}{\emph{Online}: MPC}\\[6pt]

6: & \textbf{while} control active \textbf{do} \\

7: & \quad Query $\mathbf{y}_j$, $\mathbf{u}_j$, $\mathbf{r}_j$
   \hfill $\triangleright$ Current state and reference \\

8: & \quad $\mathbf{u}_{j+1:j+m_c}
   \leftarrow \texttt{MPC}\!\left(\hat{\mathbf{F}},\mathbf{y}_j,\mathbf{u}_j,\mathbf{r}_j\right)$
   \hfill 
   \begin{tabular}[t]{@{}r@{}}
   $\triangleright$ Solve Eq.~\eqref{3.1} or Eq.~\eqref{5.1};\\
   with $\hat{\mathbf{x}}_j=\mathbf{y}_j$, $\hat{\mathbf{x}}_{j+k+1}
   =\hat{\mathbf{F}}(\hat{\mathbf{x}}_{j+k},\mathbf{u}_{j+k})$
   \end{tabular} \\

9: & \quad Observe physical state $\mathbf{y}_{j+1}$ \\

10:& \quad $j \gets j+1$ \\

11:& \textbf{end while}
\end{tabular}
}
}

\vspace{6pt}
\hrule
\end{table}





\section{Benchmark Comparison of Methods}
\label{sec:Comparison}

The proposed WSINDY--MPC framework is designed to improve robustness to noisy measurements while maintaining reliable control performance. We test it against five baseline data-driven approaches. The evaluation is conducted across varying noise levels and multiple nonlinear systems using metrics appropriate to each system, including prediction horizon accuracy, trajectory tracking error, obstacle clearance, and achieving MPC cost values. Comparisons are restricted to scenarios in which SI is feasible and a reliable MPC solution can be obtained. We compare with: (\textit{i}) SINDYc \cite{BruntonProctorKutz2016ProcNatlAcadSci}, Sec.\ref{sec:SINDYc}, (\textit{ii}) Ensemble SINDYc (E-SINDYc)\cite{FaselKutzBruntonEtAl2022ProcRSocMathPhysEngSci}, an extension of the SINDYc method that uses bootstrap aggregation (bagging) to learn the governing equations. It first applies library bagging to identify candidate terms, then uses data bootstrapping on the reduced library to obtain robust coefficient estimates, quantify their uncertainty, and enable probabilistic forecasting. We also compare with its (\textit{iii}) weak implementation (E-WSINDYc)\cite{FaselKutzBruntonEtAl2022ProcRSocMathPhysEngSci}, which uses the WSINDy method \cite{MessengerBortz2021JComputPhys}. We use the hyperparameters from the publicly available code, including $100$ bootstraps for library bagging ($90$\% sampling, $0.4$ term inclusion threshold), followed by $100$ data bootstraps for double bagging on the reduced library (with $0.6$ coefficient inclusion threshold for aggregation via median). For WSINDYc and E-WSINDYc, the parameters of the test functions are determined from data \cite{MessengerBortz2021JComputPhys,LopezNaranjoSalazarEtAl2025JComputPhys}, using the hyperparameter setting ($\tau=10^{-16}$ and $\hat{\tau}=2$) adopted in \cite{LopezNaranjoSalazarEtAl2025JComputPhys}. We also employ (\textit{iv}) dynamic mode decomposition with control (DMDc) \cite{ProctorBruntonKutz2016SIAMJApplDynSyst}, an extension of standard DMD \cite{Schmid2010JFluidMech}, that includes actuation inputs and creates low-order linear models from combined state and input data, and allows the effect of the inputs to be distinguished from the system dynamics. Finally, we utilize (\textit{v}) a multilayer NN, specifically a nonlinear autoregressive exogenous network, which models the current output as a nonlinear function of past outputs and inputs via specified delays \cite{LinHorneTinoEtAl1996IEEETransNeuralNetw}. Unless otherwise stated, we use the same NN configuration and training parameters as in \cite{KaiserKutzBrunton2018ProcRSocA}. For instance, in the NN the hyperbolic tangent sigmoid is used as the activation function, and it is first trained as a feed-forward network using the Levenberg–Marquardt algorithm (for the noise-free measurements), then the network's parameters were fixed, and it was used for prediction without any further adjustments. When there are noisy measurements, the network is trained with Bayesian regularization, which provides better robustness, but it is about four times slower than the Levenberg–Marquardt algorithm. 

For the systems with known governing equations, we also include an oracle-model MPC case, where the known model is used as the prediction model. This case provides a nominal reference to assess the performance differences that arise when the original plant model is replaced with the identified models. However, it is not expected to dominate every metric, since the nonlinear MPC problem is generally nonconvex and may converge to local solutions \cite{MayneRawlingsRaoEtAl2000Automatica}.

For consistency with the SINDYc-based methods, all available data points are used as query points in the WSINDYc-based methods. Note that, for both SINDy and E-SINDy, the time derivatives are approximated using a centered fourth-order finite difference \cite{KaiserKutzBrunton2018ProcRSocA}.

\subsection{Plasma boundary control}
High-fidelity plasma simulations are computationally expensive for real-time control. Moreover, plasma boundary data is noisy, motivating the use of robust data-driven identification methods to obtain reliable reduced-order models. In this sense, we first provide a comparison using data generated from the SOLPS-ITER (Scrape‑Off Layer Plasma Simulation)
 code \cite{BonninDekeyserPittsEtAl2016PlasmaandFusionResearch}. The configuration is shown in  Fig.~\ref{fig_SOLPS_mesh}, and a brief description of the governing equations is presented in Appendix \ref{appendix:plasma}.

\begin{figure}[tb]
\centering\includegraphics[width=2.5in]{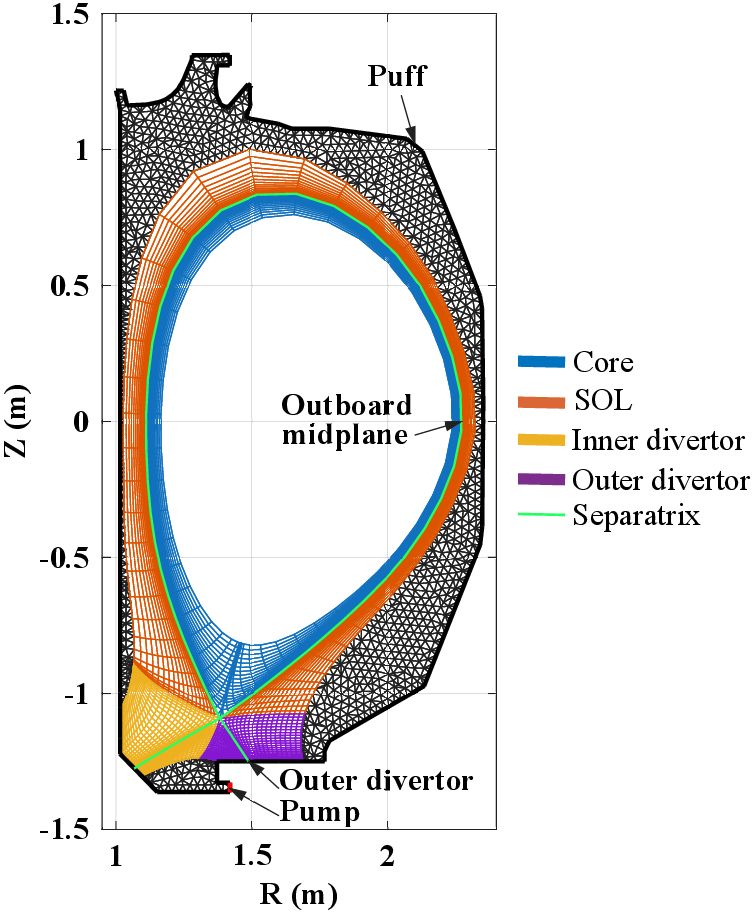}
\caption{SOLPS-ITER configuration and mesh, adapted from \cite{LoreDePascualeLaiuEtAl2023NuclFusion}, CC BY 4.0. Colored regions denote various plasma regions. (Online version in color.)}
\label{fig_SOLPS_mesh}
\vspace{-8mm}
\end{figure}

\begin{figure}[tb]
\centering\includegraphics[width=5.0in]{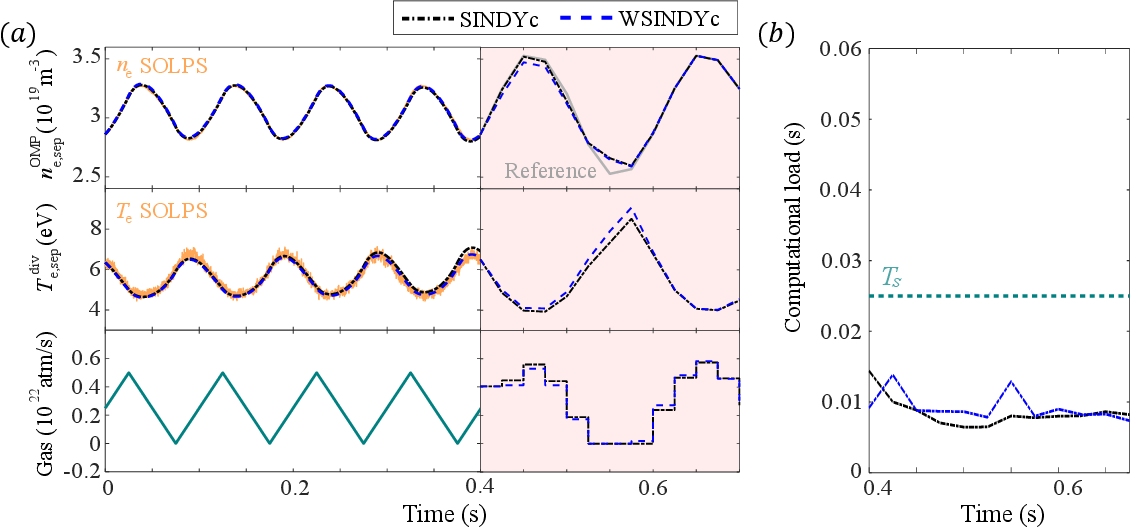}
\caption{($a$) Time series and models identified for the SOLPS-ITER are presented in the white region, while the MPC stage is represented in the pink region. The top rows are the states of interest, while the bottom row is the control input, ($b$) per-MPC-step optimization time contrasted with actuator interval $T_s$. (Online version in color.)}
\label{fig_SOLPS_TS}
\vspace{-7mm}
\end{figure}

The time series $n^\text{{OMP}}_\text{e,sep}$ (electron density at the outer midplane separatrix) and $T^{\text{div}}_\text{e,sep}$ (electron temperature at the divertor separatrix), shown in Fig.~\ref{fig_SOLPS_TS}a as continuous lines,  have noise levels of approximately 0.1$\%$ and 15$\%$, respectively. This noise is primarily due to Monte Carlo noise introduced by the EIRENE simulation. To identify the models that describe the plasma boundary response to gas-puff actuation, SINDYc-based methods are constructed using linear models \cite{LoreDePascualeLaiuEtAl2023NuclFusion}. As can be seen in Fig.~\ref{fig_SOLPS_TS}a, both SINDYc and WSINDYc provide a similar and accurate representation of the dynamics of the two states. We then carried out MPC  using the models provided by SINDYc and WSINDYc. Figure~\ref{fig_SOLPS_TS}a (pink region) shows that the reference state $x_1$ is well traced for both methods. The divertor temperature is not explicitly tracked and evolves according to the coupled dynamics. Different from the implementation of \cite{KaiserKutzBrunton2018ProcRSocA}, where each MPC prediction step advances the model by a single numerical integration step, we integrate the model over each MPC update interval $T_s$ using internal integration steps $n_s^\text{M}=T_s/\Delta t^{\text{M}}$, while keeping the control input $\mathbf{u}$ constant over the interval. This preserves the same physical control and prediction horizons $T_c=T_p=m_cT_s$, but improves the numerical accuracy of the predicted state evolution within each MPC step. For this system, the resulting finer temporal resolution provides a better representation of $n^\text{{OMP}}_\text{e,sep}$ and $T^{\text{div}}_\text{e,sep}$, which may not be adequately captured using a single integration step per prediction interval.

In real-time control, for safety and appropriate performance, each per-MPC-step optimization time must finish below the actuator interval $T_s$. Figure~\ref{fig_SOLPS_TS}b presents the computation times for both methods, which remain consistently low throughout. Note that to approach real-time requirements, we compiled the MPC objective from \cite{KaiserKutzBrunton2018ProcRSocA} to C via MATLAB Coder, reducing the per-MPC-step solve time to the order of $T_s$.

In real-world applications, measurements are almost always contaminated by noise, which can significantly impact the system identification accuracy and control performance. In this particular simulated system, the noise arises primarily from the Monte Carlo neutral transport solver that is part of SOLPS-ITER. Reducing the noise requires increasing the number of simulated particles, reflecting the inherent trade-off between simulation fidelity and computational cost \cite{LoreDePascualeLaiuEtAl2023NuclFusion}. To assess robustness, increasing levels of synthetic noise are added to the noisy states data of Fig.~\ref{fig_SOLPS_TS}a, before SI, with noise magnitudes  $\eta_i \in (0.01,1)$.

\begin{figure}[tb]
\centering\includegraphics[width=5.25in]{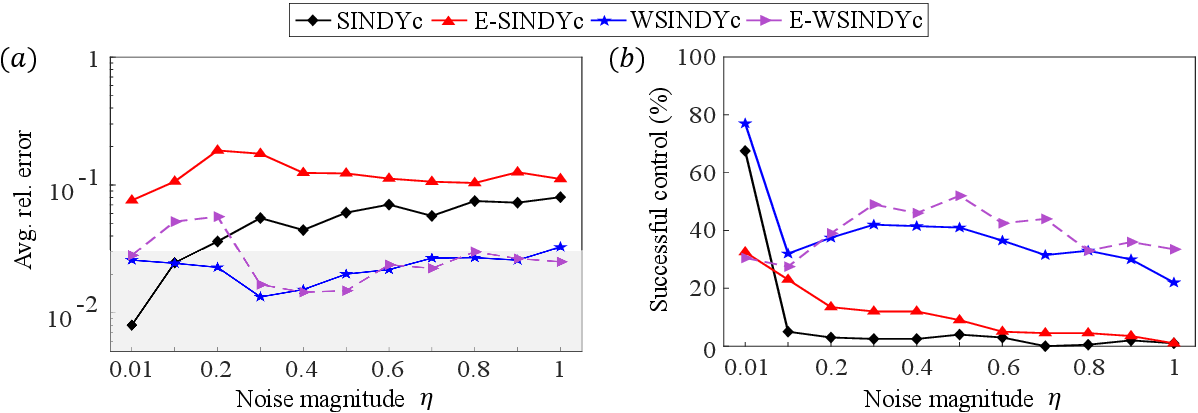}
\caption{Performance for different noise levels for the SOLPS-ITER example: ($a$) Average relative error between the obtained state and the reference signal, ($b$) Success rate. E-WSINDYc demonstrates better tracking of state $x_1$ with generally lower relative error and higher success rates. Beyond the shaded region, performance degrades: models diverge and have no predictive power. (Online version in color.)}
\label{fig_SOLPS_noise_dep}
\vspace{-5mm}
\end{figure}

For each noise level, reduced-order models are identified using a library of candidate functions constructed with polynomial terms up to second order, rather than the linear model of the previous demonstration, to capture the nonlinear plasma response, and subsequently used within the MPC framework. Control performance is evaluated by comparing the tracked first state against its reference signal using the average relative error. In addition, a success rate is computed as the fraction of realizations achieving an average relative error below $3\%$, chosen as a representative tolerance for acceptable tracking performance. To improve readability, Fig.~\ref{fig_SOLPS_noise_dep} displays the medians, while the bars representing the interquartile range, from the first to the third quartile, are depicted in Fig.~\ref{fig_SOLPS_noise_dep_box} in Appendix~\ref{appendix:plasma}. As shown in Figure~\ref{fig_SOLPS_noise_dep}, E-WSINDYc provides the best result for MPC, followed by WSINDYc, E-SINDYc, and finally SINDYc. This not only confirms the repeated success of the weak-form methods that have recently been developed to tackle the presence of noise \cite{MessengerBortz2021MultiscaleModelSimul,BortzMessengerDukic2023BullMathBiol,MessengerBortz2022PhysicaD,BortzMessengerTran2024NumericalAnalysisMeetsMachineLearning,MessengerBortz2021JComputPhys,LopezNaranjoSalazarEtAl2025JComputPhys,TranBortz2025arXiv250703206}, but also that the ensemble method \cite{FaselKutzBruntonEtAl2022ProcRSocMathPhysEngSci} combined with WSINDYc is a robust approach for MPC. The per-MPC-step optimization times of Eq.~\eqref{3.1} are presented in Fig.\ref{fig_SOLPS_MPC_time} in Appendix \ref{appendix:plasma}.

\subsection{Tracking a drone}
\begin{figure}[tb]
\centering\includegraphics[width=5.0in]{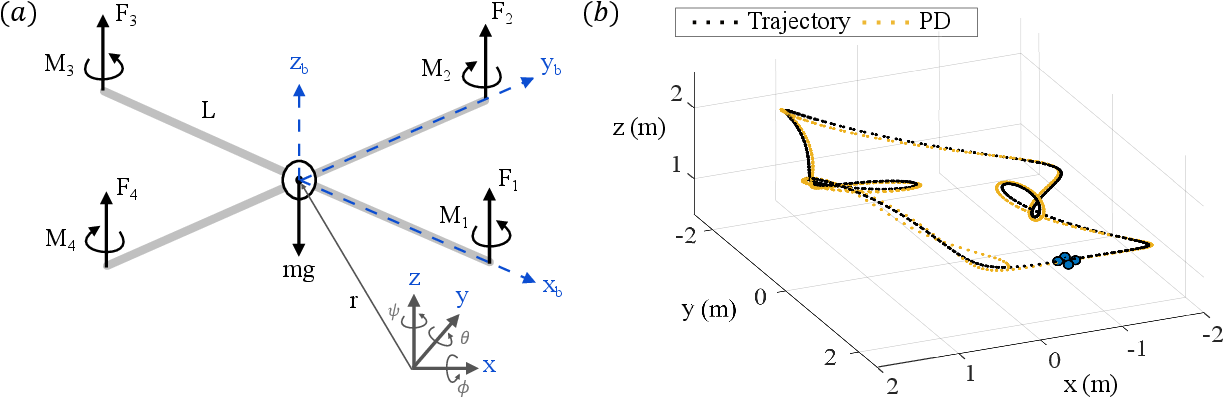}
\caption{($a$) Inertial (solid lines), body (dashed lines) reference frames, along with forces and moments acting on quadrotor, ($b$) Data generation using the PD controller (Online version in color.)}
\label{fig_drone_training_traj}
\vspace{-7mm}
\end{figure}

In this subsection, we consider a quadrotor simulator \cite{MichaelMellingerLindseyEtAl2010IEEERobotAutomatMag} tasked with tracking trajectories in the presence of a spherical obstacle \cite{LeeKimKimEtAl2025IETControlTheoryAppl}. In this system, the framework consists of three stages: an offline stage, where a proportional-derivative (PD) controller is used to track reference trajectories chosen to sufficiently excite the quadrotor simulator, thus generating the state and input data required for system identification. This controller is not used for control design in this work; it serves solely as a data generation mechanism. In the subsequent stage, the collected state and effective input data are used to identify a model of the quadrotor simulator using the different identification methods considered. In the final stage, the identified models are used within the MPC configuration to perform trajectory tracking with obstacle avoidance. Figure~\ref{fig_drone_training_traj}a shows the configuration of the vehicle's dynamic model, and Fig.~\ref{fig_drone_training_traj}b presents the reference trajectory and the actual flight path achieved by the PD controller. The system equations and simulation parameters are given in Appendix \ref{appendix:drone}.

To perform SI of the quadrotor simulator, we follow the same approach as in \cite{LeeKimKimEtAl2025IETControlTheoryAppl}, which uses a SINDy-based method. The details of SI are presented in Appendix \ref{appendix:drone}. Regarding NN being applied in this system, it exhibits limited system identification accuracy and lacks the robustness required for reliable MPC. A similar observation was found in a similar system in \cite{LiuHongPiprekEtAl2025IEEETransAerospElectronSyst}. 

To carry out MPC for tracking trajectory in the presence of an obstacle, we modified Eq.~\eqref{3.1} as
\begin{align}\label{5.1}
J &= 
\sum_{k=0}^{m_p-1} 
\Bigg(\left\| \hat{\mathbf{x}}_{j+k} - \mathbf{r}_k \right\|_{\mathbf{Q}}^{2}
    + Q_{\text{obs}}
    \Big[\max\!\big(0,\; D_{\min} - \| \hat{\mathbf{p}}_{j+k} - \mathbf{o} \|_2 \big)
    \Big]^{2}
\Bigg) \nonumber \\
&\quad +\sum_{k=1}^{m_c-1}\Big(\left\| \mathbf{u}_{j+k} \right\|_{\mathbf{R_u}}^{2}
    + \left\| \Delta\mathbf{u}_{j+k} \right\|_{\mathbf{R}_{\Delta u}}^{2}\Big).
\end{align}
where $\hat{\mathbf{p}}_{j+k}\in\mathbb{R}^{3}$ is the predicted position, $\mathbf{o}\in\mathbb{R}^{3}$ is the center of the obstacle, $D_\text{min}$ is a safety radius set to at least the sum of the obstacle radius and the drone's arm length, defining the minimum clearance required for collision-free flight, and $Q_\text{obs}$ is a penalty weight, which is set large enough to enforce avoidance. As in the previous system, during the optimization, the dynamics are integrated over a full control interval per MPC step using multiple smaller integration substeps, allowing the prediction horizon to more accurately capture the nonlinear flight behavior of the drone.

\begin{figure}[tb]
\centering\includegraphics[width=5.0in]{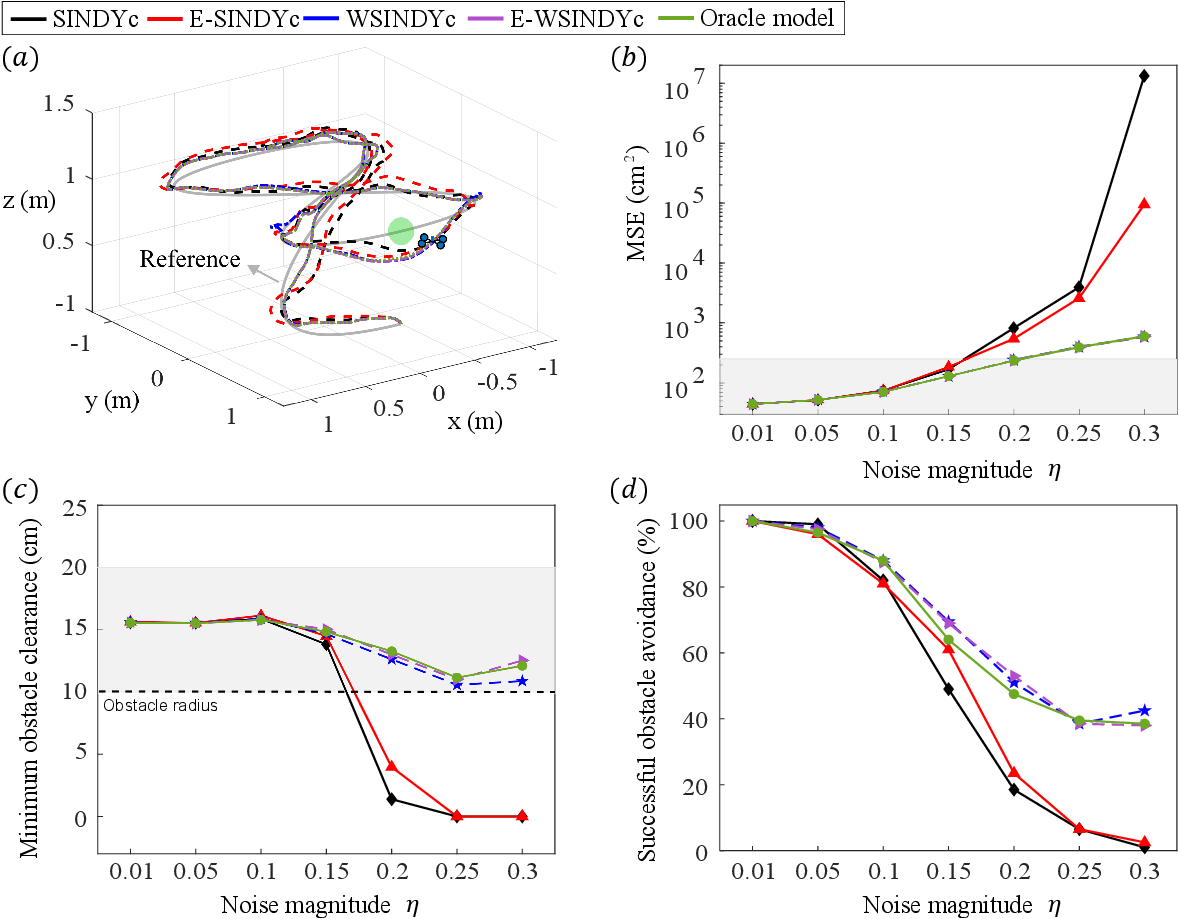}
\caption{($a$) Trajectory tracking performance and obstacle avoidance, $\eta=0.15$ and the obstacle radius is $10~\text{cm}$; Performance for different noise levels: ($b$) mean square error, and ($c$) minimum obstacle clearance, for 200 noise realizations each. Outside the shaded region, prediction accuracy deteriorates, and the resulting control actions fail to guarantee safe obstacle avoidance.  (Online version in color.)}
\label{fig_drone_dep}
\vspace{-7mm}
\end{figure}

Figure~\ref{fig_drone_dep}a illustrates representative trajectories for a single realization at a fixed noise level, including the reference trajectory (solid grey), the spherical obstacle (green), and the trajectories generated by MPC using the different identified models. WSINDYc and E-WSINDYc provide the best and similar results, followed by E-SINDYc, and finally SINDYc. To evaluate robustness across noise levels, the tracking error is quantified using MSE computed only outside the obstacle-influenced region, defined by the reference points whose distance to the obstacle center exceeds $D_\text{min}$. The reference trajectory is adapted from \cite{Lopez-SanchezMoreno-Valenzuela2023AnnualReviewsinControl}. Figure~\ref{fig_drone_dep}b shows that WSINDY--MPC and E-WSINDY--MPC outperform the other two methods in high noise levels. The corresponding bars representing the interquartile range, from the first quartile to the third quartile, for Figure\ref{fig_drone_dep}b are presented in figure \ref{fig_drone_dep_box}a in Appendix \ref{appendix:drone}.

Although tracking error captures trajectory accuracy, it does not reflect collision risk. Therefore, we compute the minimum effective clearance between the obstacle surface and the closest propeller tip along each trajectory, accounting for the arm length $L$ when calculating the three-dimensional distance to the obstacle. Figure~\ref{fig_drone_dep}c depicts that up to $\eta=0.15$, the medians remain above the obstacle radius, but the spread shown in Fig~\ref{fig_drone_dep_box}b in Appendix \ref{appendix:drone} indicates that some realizations of SINDYc and E-SINDYc lose safe clearance. From  $\eta=0.2$ onward, the medians of SINDYc and E-SINDYc drop, indicating frequent collisions, whereas WSINDYc and E-WSINDYc maintain larger safety margins throughout the tested range. This is corroborated in Fig.~\ref {fig_drone_dep}d, which presents the number of flights that successfully avoided the obstacle, defined as those whose minimum clearance falls within the gray region of $10$ to $20$ cm. Figure~\ref{fig_drone_MPC_time} in Appendix \ref{appendix:drone} displays the per-MPC-step solve times of Eq.~\eqref{5.1}.

\subsection{Lorenz attractor}

 We now compare performances using the Lorenz $63$ system, a classic example of deterministic chaos that arises from a simplified model of atmospheric convection\cite{Lorenz1963JAtmosSci}. The governing equations and simulation parameters are presented in Appendix \ref{appendix:lorenz}.

\begin{figure}[tb]
\centering\includegraphics[width=5.0in]{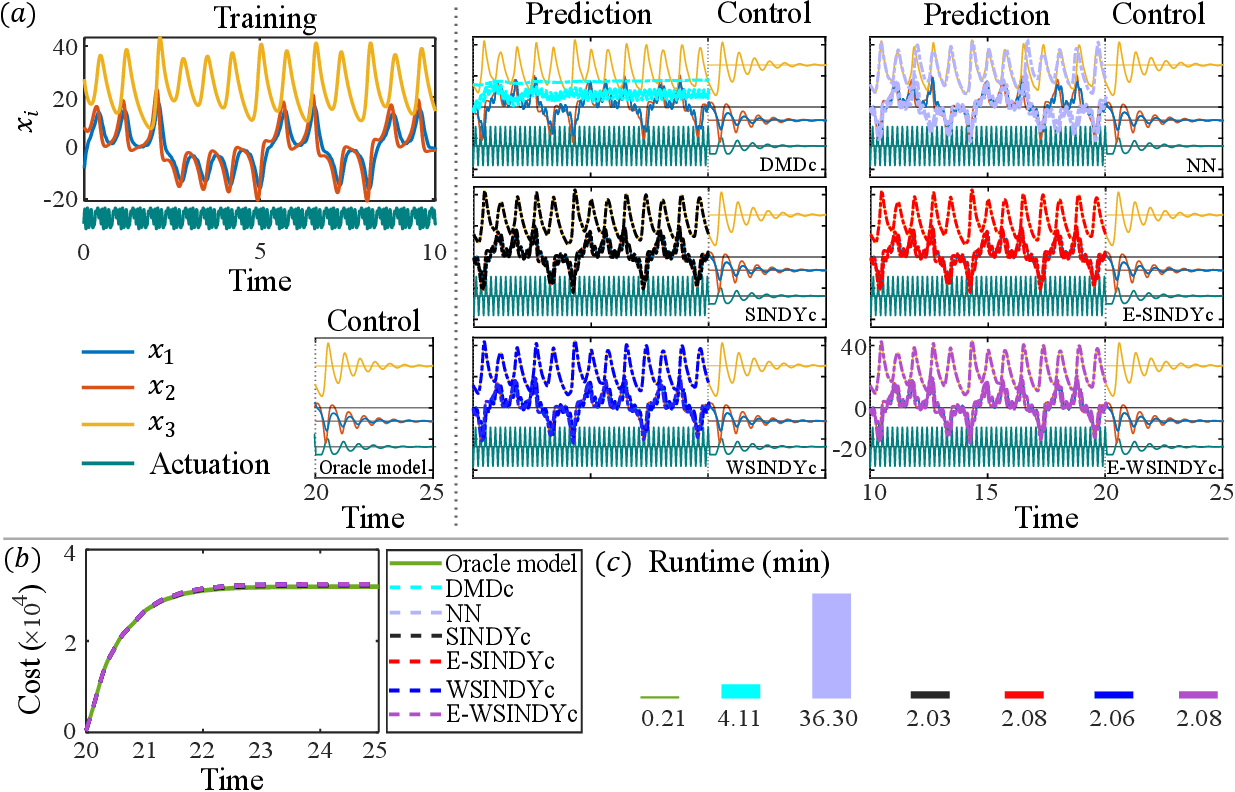}
\caption{Prediction and control results for the Lorenz attractor: ($a$) time series of states and inputs during training, validation, and control stage, ($b$) Cumulative cost of Eq.~\eqref{3.1}, and ($c$) Runtime of the MPC optimization process. (Online version in color.)}
\label{fig_Lorenz_1}
\vspace{-5mm}
\end{figure}

Figure~\ref{fig_Lorenz_1}a (training) shows the simulated time series. The models were trained over the first $10$ time units (10,001 discrete points at time step $10^{-3}$). Ranking model complexity, DMDc provides the least complex model, as it assumes a linear mapping between consecutive system states, capturing dominant modes of dynamics through linear superposition. SINDYc, E-SINDYc, WSINDYc, and E-WSINDYc, are intermediate, due to the interpretable equations as a result of the sparse selection of nonlinear terms from the candidate functions of a library. NN is the least interpretable, since it acts as a black box whose trainable parameters do not correspond to physically meaningful terms. 

After training, the models are validated over the subsequent 10 time units with the sinusoidal forcing $u(t)=(5\sin(30t))^3$, where Figure~\ref{fig_Lorenz_1}a (prediction) shows that the SINDYc-based methods provide good results, followed by the NN, and finally DMDc. Subsequently, MPC is implemented for $5$ time units. Figure~\ref{fig_Lorenz_1}a (control) displays that the control performed by all models presents similar results. Regarding the MPC calculation time (Figure~\ref{fig_Lorenz_1}c), NN needs the longest time, followed by DMDc, and last the SINDYc-based methods. The longer MPC computation time for the NN is mainly due to the repeated closed-loop evaluation of the model at each optimization iteration. In our MATLAB implementation, this sequential evaluation introduces substantially more overhead than the SINDYc-based models, which require only polynomial evaluations, or DMDc, which involves a linear matrix update. The oracle-model MPC provides the fastest time, as expected, because the exact Lorenz equations are cheap to evaluate. Note that with the C implementation of the MPC objective functions, the MPC optimization runtime with the oracle, DMDc, and sparse-regression models is approximately $12$ seconds. We report the time for NN only in the native form, since MATLAB Coder requires an explicit closed-form predictor, which a trained network does not provide.

\begin{figure}[tb]
\centering\includegraphics[width=5.0in]{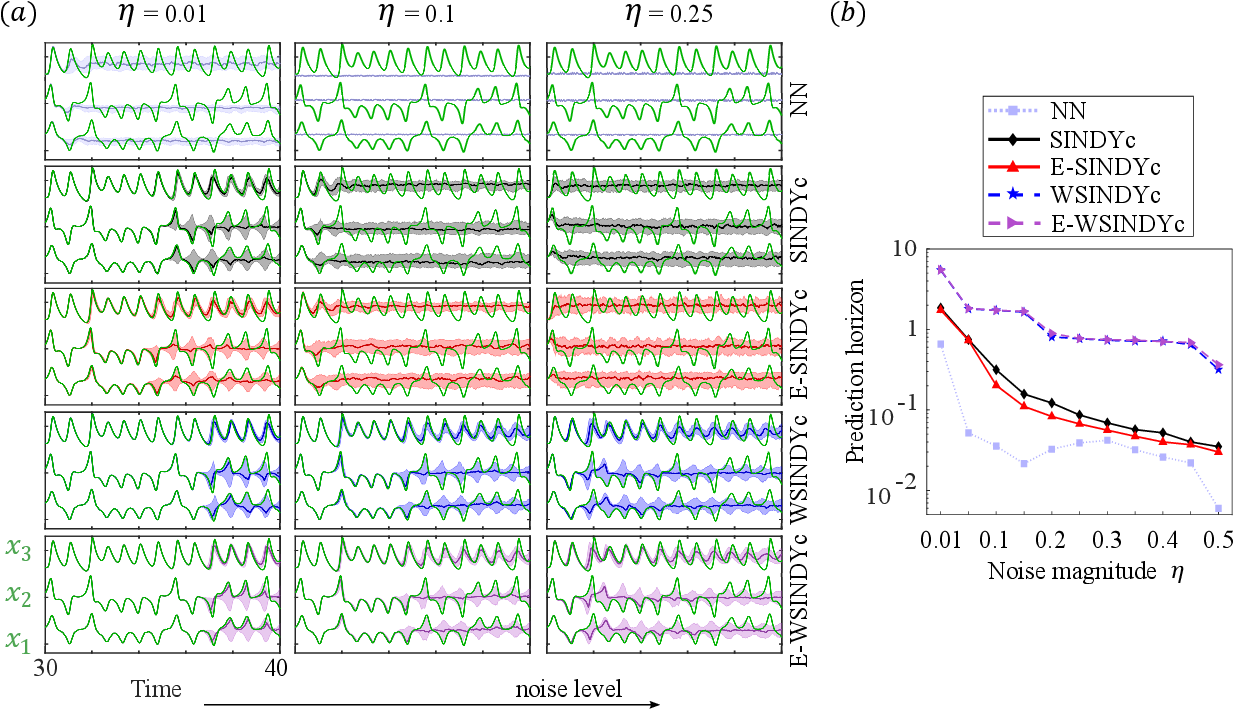}
\caption{Validation prediction performance for different noise levels: ($a$) time series, showing median (thick colored line) and 25–75th percentiles (shaded region), and ($b$) prediction horizon in time units, from 200 noise realizations each (Online version in color.)}
\label{fig_Lorenz2}
\vspace{-7mm}
\end{figure}

We evaluate the prediction performance across $200$ independent noise realizations. Figure~\ref{fig_Lorenz2}a shows the validation forecasting for NN and the SINDYc-based methods using noise levels $\eta = 0.01, 0.1,$ and $0.25$. As the noise level increases, all the models' performance degrades, though WSINDYc consistently outperforms the other three methods. To quantify the prediction performance, the prediction horizon metric \cite{KaiserKutzBrunton2018ProcRSocA} is used, which is the point in time when the prediction error exceeds a chosen tolerance, $\sqrt{\sum_{i=1}^{3}(x_i - \hat{x}_i)^2} < \varepsilon$, with the tolerance set to $\varepsilon=3$, corresponding to roughly $10\%$ of the magnitude of each state variable (given states on the order of $O(10^1)$). Beyond this point, the predicted states visibly diverge from the true ones. Figure~\ref{fig_Lorenz2}b depicts the superiority of WSINDYc-based methods to predict the states, with E-WSINDYc slightly providing better results than WSINDYc. The interquartile ranges for Figures \ref{fig_Lorenz2}b are shown figure \ref{fig_Lorenz2_box}a in Appendix \ref{appendix:lorenz}.

\begin{figure}[tb]
\centering\includegraphics[width=5.0in]{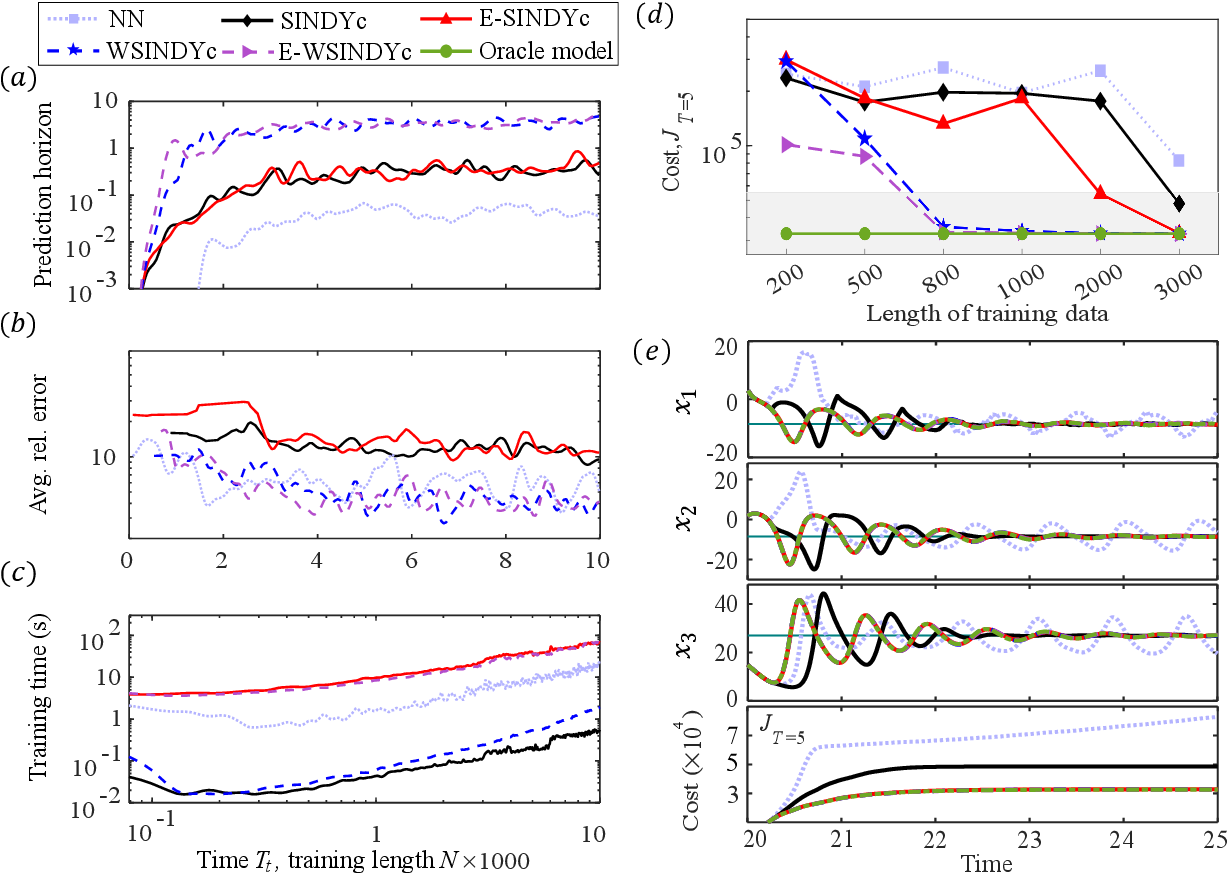}
\caption{Validation prediction and control performance on the Lorenz system as the length of training data increases, when $\eta = 0.1$: ($a$) prediction horizon, ($b$) average relative prediction error, ($c$) training time, ($d$) terminal cumulative cost, using Eq.~\eqref{3.1}, over $5$ time units, where from $800$ data points onward, E-WSINDYc demonstrates more consistent control performance with generally lower costs, and ($e$) time series of states and cost of the best model at $3000$ data points for all methods. Beyond the shaded region, performance degrades, with some models even diverging and most lacking predictive power. (Online version in color.)}
\label{fig_Lorenz3}
\vspace{-7mm}
\end{figure}

In SI and MPC, having sufficient data is important, as they determines the accuracy of the identified model and the effectiveness of control. We now examine the performance of the WSINDYc method in different amounts of data when $\eta=0.1$. As the data length increases,  Fig.~\ref{fig_Lorenz3}a-c shows the corresponding prediction horizon, average relative error, and training time. When training data are limited, WSINDYc and E-WSINDYc demonstrate superior performance compared to NN, SINDYc, and E-SINDYc, achieving a prediction horizon of approximately $3$ time units and a lower average relative error. With large datasets, WSINDYc and E-WSINDYc perform similarly. 

To assess the computational cost of each method, we measured the training time $\tau$ as a function of $N$ and performed a power-law fit $\tau \sim N^\alpha
$ in a log-log space (using MATLAB's \texttt{polyfit} function, and shown in Fig.~\ref{fig_Lorenz3}c). The results indicate the SINDYc and NARX exhibit approximately linear scaling ($\alpha=0.98$ and $\alpha=0.92$, respectively), while the ensemble variants E-SINDYc and E-WSINDYc scale sublinearly ($\alpha=0.78$ and $\alpha=0.82$). In contrast, WSINDYc displays a superlinear scaling ($\alpha=1.28$). To enhance readability, the plotted curves were downsampled by a factor of $5$; outliers are removed using MATLAB's \texttt{rmoutliers} function, and the remaining data are smoothed with \texttt{smoothdata} using a Gaussian window of length $10$. 

We next evaluate the models in MPC, where a new control input is calculated every $10$ time steps and remains fixed until the next update. For Fig.~\ref{fig_Lorenz3}d, $200$ distinct sets of noise data are created as input for training, and the resulting identified systems are used for the evaluation of the MPC of the five methods. 
From 800 data points onward, both E-WSINDYc and WSINDYc obtain $J$ values within the shaded region, indicating effective control. We compute the median $J$ for each model, and select the training length with the minimum median value. The per-MPC-step optimization times of Eq.~\eqref{3.1} are presented in Fig. ~\ref{fig_Lorenz4} in Appendix \ref{appendix:lorenz}. Fig.~\ref{fig_Lorenz3}e then shows the corresponding state trajectories and terminal cumulative cost over $5$ time units for those selected cases. Overall, WSINDYc and E-WSINDYc show superior performance on large datasets, followed by E-SINDYc, SINDYc, and last NN; with E-WSINDYc achieving slightly better results than WSINDYc. The interquartile ranges for Figs.~\ref{fig_Lorenz3}d are shown in Fig.~\ref{fig_Lorenz2_box}b in Appendix \ref{appendix:lorenz}.

\maketitle

\section{Conclusion}
\label{sec:Conc}

In this paper, we presented WSINDY--MPC, a data-driven MPC framework that exhibits strong robustness in the presence of high noise levels. The method is demonstrated on four case studies: a fusion model, a quadrotor model, the chaotic Lorenz $63$ system, and a non-affine F8 Crusador model (Appendix \ref{appendix:f8}). Across these examples, both WSINDY--MPC and E-WSINDY--MPC frameworks can control these strongly nonlinear systems using only measurement data  and achieve similar levels of performance. While E-WSINDY--MPC provides modest gains in certain scenarios, WSINDY--MPC outperforms it in others and trains faster, indicating that their relative performance is problem dependent. Overall, the weak formulation plays a central role in achieving improved performance across most metrics, enabling both approaches to outperform alternative data-driven MPC methods, including DMDc, NN, SINDYc, and E-SINDYc.

While SINDYc can operate in the low-data limit by exploiting sparsity in the identified model, WSINDYc relies on test functions supported on multiple consecutive time points \cite{MessengerBortz2021MultiscaleModelSimul}. This introduces a minimum data requirement that can exceed that of derivative-based approaches, particularly for very short trajectories. In that regard, advanced frameworks to construct test functions will be worth exploring for extending WSINDYc to data-limited regimes\cite{TranBortz2025arXiv250703206}.

Further investigation of the WSINDY--MPC framework will focus on applying it to experimental data, building on recent work that identifies opportunities for sparse regression in plasma physics for model discovery for tokamak plasma boundary \cite{KaptanogluHansenLoreEtAl2023PhysPlasmas}. Also, it will be interesting to apply the proposed methodology to drone data \cite{GuevaraVarela-AldasGandolfoEtAl2025Drones,OsmanXiaMahdiEtAl2025IntlJRobustNonlinear}, where the presence of noise and abrupt changes during operation are important aspects to consider. 
Beyond the single-drone setting considered here, extending our framework to drone swarms is a potential research area \cite{YazdannikTayefiFarrokh2023arXiv231112906},
with recent work exploring communication-efficient online path control\cite{ShiriParkBennis2020IEEETransCommun}. Since real-world applications demand the controller to run in real time, we compiled the MPC objective and constraints functions into C, bringing the per-step solve times to the order of $T_s$. Achieving the margin for reliable real-time operation will require further acceleration, for which a MPC framework such as \texttt{acados}\cite{FreyBaumgartnerFrisonEtAl2025OptimControlApplMethods} is a natural direction for future work.

In all cases considered here, the discrete-time models are identified assuming that the actuation input is known. Extending the proposed methodology to handle unknown inputs, as demonstrated in \cite{ZhangGuoZhangEtAl2024PhilTransRSocA}, represents a promising direction for future research. Moreover, the identified models lead to ODEs, which are subsequently used for MPC. Incorporating systems governed by partial differential equations into this framework, and integrating them with MPC, is another important direction \cite{KordaMezic2018Automatica}. Additionally, the synthetic noise considered is Gaussian; it would be of interest to investigate the performance of the framework under alternative noise distributions \cite{LopezNaranjoSalazarEtAl2025JComputPhys}. Furthermore, incorporating online system identification would enable continuous model refinement during closed-loop operation while accounting for parametric uncertainties \cite{ManzoorPeiSunEtAl2022Drones}. Addressing systems with nonlinear-in-parameters ODEs constitutes another relevant direction for future work \cite{RummelMessengerBeckerEtAl2025arXiv250208881}.

Finally, while our work demonstrates an empirical closed-loop performance, developing robust MPC formulations with rigorous stability analysis \cite{WuTranRinconEtAl2019AIChEJournal,WuTranRinconEtAl2019AIChEJournala,WuChristofidesWuEtAl2025RevChemEng,AbdullahWuChristofides2021ComputChemEng,AbdullahChristofides2023ChemicalEngineeringResearchandDesign} is a valuable path to pursue. Such extensions could provide guarantees of recursive feasibility, state boundedness, and modeling error, which would further strengthen the framework for real-world deployment.

\vskip6pt

\enlargethispage{20pt}

\dataccess{The code used in this work is available at: \url{https://github.com/MathBioCU/WSINDY-MPC} and can be used to regenerate all figures in this work.}

\aucontribute{All authors conceived of the work, designed the study, and edited the manuscript. C.L. carried out the computations and wrote the first draft of the manuscript.}
\competing{We declare we have no competing interests.}
\funding{C.L., A.J.C., S.B., and D.M.B. acknowledge support from the Department of Energy Mathematical Multifaceted Integrated Capability Center (MMICC) grant to Michigan State University (DE-SC0023164) and to the University of Colorado (DE-SC0023346). A.J.C. would also like to thank DOE/NNSA grant DE-NA0004265. This manuscript has been authored in part by UT-Battelle, LLC, under contract DE-AC05-00OR22725 with the US Department of Energy.  This work also utilized the Blanca condo computing resource at the University of Colorado Boulder. Blanca is jointly funded by computing users and the University of Colorado Boulder.}

\appendix

\section*{Appendix A. Plasma boundary dynamics and simulation parameters}
\refstepcounter{section}\label{appendix:plasma}
\setcounter{equation}{0}
For this system, coupled two-dimensional plasma fluid equations are solved for the tokamak boundary region. A detailed description of the SOLPS-ITER equations can be found in \cite{BonninDekeyserPittsEtAl2016PlasmaandFusionResearch} and references therein. Here, following Ref. \cite{LoreDePascualeLaiuEtAl2023NuclFusion}, we present a simplified version of the relevant equations.

\begin{align}\label{4.1}
\frac{\partial n}{\partial t} + \nabla \cdot \Gamma_n &= S_n, \tag{A 1a}\\
\mathrm{m} \frac{\partial (n V_{\parallel})}{\partial t} + \nabla \cdot \Gamma_{\mathrm{m}} &= -\nabla \mathrm{p} + \mathrm{q}nE_{\parallel} + S_{\mathrm{m}}, \tag{A 1b}\\
\frac{3}{2} \frac{\partial (nT)}{\partial t} + \nabla \cdot \left( \mathbf{q} + \frac{3}{2} \mathrm{p} \mathbf{V} \right) &= -\mathrm{p} \nabla \cdot \mathbf{V} + \mathcal{C} + S_E\tag{A 1c}.
\end{align}

\noindent where $n$ is the plasma particle density, $T$ is the plasma temperature, $V$ is the plasma flow velocity; $\Gamma_n$ and $\Gamma_{\mathrm{m}}$ are the particle and parallel momentum flux densities, respectively; $\mathbf{q}+\frac{3}{2} \mathrm{p} \mathbf{V}$ the internal energy flux density, with $\mathbf{q}$ the heat flux and $\mathrm{p}=nT$ the static pressure. The term $\mathcal{C}$ represents collisional energy exchange between the plasma species, and $S$ denotes the source terms arising from plasma-neutral interactions, which are computed by EIRENE \cite{ReiterBaelmansBorner2005FusionScienceandTechnology}, a kinetic neutral solver that uses a Monte Carlo approach on an unstructured mesh to model plasma–neutral and plasma–material interactions (see Fig.~\ref{fig_SOLPS_mesh}).

To obtain the training data, two state variables are extracted from the full two-dimensional time series data of the SOLPS-ITER simulations: the separatrix electron density at the outer midplane, $n^\text{{OMP}}_\text{e,sep}$, and the separatrix electron temperature at the outer divertor, $T^{\text{div}}_\text{e,sep}$ \cite{LoreDePascualeLaiuEtAl2023NuclFusion}. The density $n^\text{{OMP}}_\text{e,sep}$ is an important operational quantity for tokamak control. It is directly actuated by main ion gas puffing $u$ (deuterium injection) and is strongly linked to the overall divertor conditions. The $T^{\text{div}}_\text{e,sep}$ serves as a crucial indicator for the divertor integrity and scales proportionally with a large number of quantities that characterize divertor heat and particle fluxes. The time series of plasma boundary data was generated under time-varying $u$ with a duration of 0.4 s and was uniformly sampled, resulting in 4001 data points. As indicated in \cite{LoreDePascualeLaiuEtAl2023NuclFusion}, the input $u$ and state variables $x_1=n^\text{{OMP}}_\text{e,sep}$ and the inverse divertor electron temperature $x_2=1/T^{\text{div}}_\text{e,sep}$ were normalized  by $10^{18}$, $10^{16}$, and $10^{-4}$, respectively, in order to improve the accuracy of the models.

We perform MPC over a 0.296 s interval. The gas-puff actuator is updated in constant piecewise steps of $T_s=0.025$ s, a new control input is calculated every $250$ time steps and remains fixed until the next update, and the control objective is defined by a sinusoidal reference trajectory for $n^\text{{OMP}}_\text{e,sep}$. We define its initial value as the last value of state $x_1$ used for training, and it oscillates at a frequency of 5 Hz.  We set $m_p=12$ following \cite{KaiserKutzBrunton2018ProcRSocA} and $m_c=4$. The model timestep $\Delta t^\text{M}$ and the system's $\Delta t^\text{plant}$ equal to $10^{-4}$ s, matching the SOLPS-ITER simulation timestep \cite{LoreDePascualeLaiuEtAl2023NuclFusion}, and $n_s^{\text{plant}}=1$. Only $n^\text{{OMP}}_\text{e,sep}$ is tracked ($\mathbf{Q} = \text{diag}(25,0)$), as the two state variables are coupled \cite{LoreDePascualeLaiuEtAl2023NuclFusion}. Small penalties $\mathbf{R}_u = \mathbf{R}_{\Delta u} = \text{diag}(0.001,0.001)$ prioritize tracking over actuator effort, and we use the constraint $u \in [0, 6\times10^3]$ to enforce physical nonnegativity of the gas puff rate while allowing sufficient actuation range. Both DMDc and NN were excluded from this system, as they did not accurately capture the system dynamics during SI, precluding its use for MPC.

\begin{figure}[tb]
\centering\includegraphics[width=5.0in]{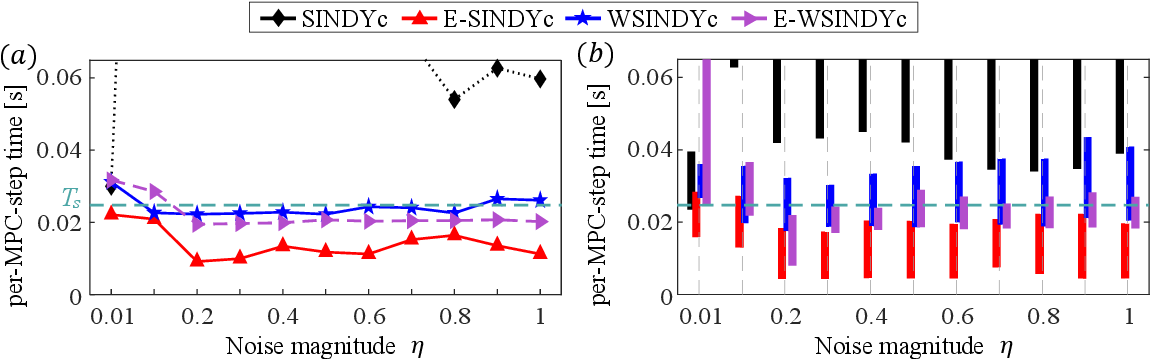}
\caption{SOLPS, per-MPC-step runtime using Eq.~\eqref{3.1}: ($a$) Medians ($b$) first to third quartile of ($a$), from 200 noise realizations each (Online version in color.)}
\label{fig_SOLPS_MPC_time}
\vspace{-7mm}
\end{figure}

Figure \ref{fig_SOLPS_MPC_time} displays the per-MPC-step solve time. Only E-SINDY--MPC guarantees that computation terminates on schedule, while the weak-form-based methods provide the solve times close to $T_s$.

\begin{figure}[tb]
\centering\includegraphics[width=3.25in]{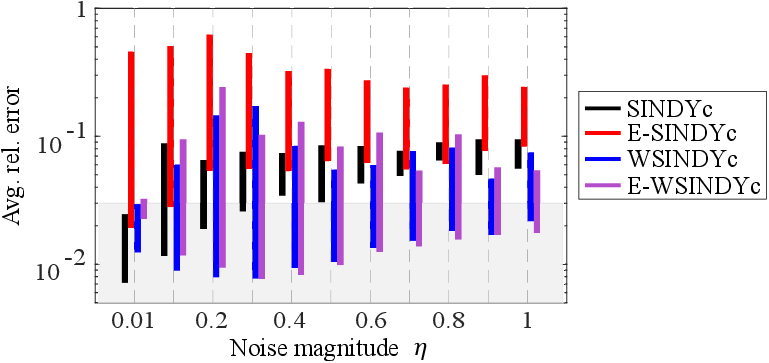}
\caption{Performance for different noise levels: MSE between the obtained state and the reference signal, from first to third quartile. E-WSINDYc generally shows less dispersion in the data. Beyond the shaded region, performance degrades, some models diverge, and have no predictive power. (Online version in color.)}
\label{fig_SOLPS_noise_dep_box}
\vspace{-7mm}
\end{figure}

\section*{Appendix B. Drone dynamics and simulation parameters} \refstepcounter{section}\label{appendix:drone}
\renewcommand{\theequation}{B \arabic{equation}}
\setcounter{equation}{0}
The acceleration of the center of mass is governed by the equation

\begin{align}\
m\mathbf{\ddot{p}} = m\mathbf{g} + \mathbf{R}(\boldsymbol{\beta})
\begin{bmatrix} 
0 \\ 
0 \\ 
F 
\end{bmatrix},
\end{align}
where $m$ is the mass of the quadrotor, $\mathbf{p}=[x\ y\ z]^\top\in\mathbb{R}^3$ is the inertial position, $\mathbf{g}=[0,0,-g]^{\top}$, with $g$ the acceleration due to gravity, $\mathbf{R}(\boldsymbol{\beta})$ the rotation matrix for transforming the coordinates from the body to inertial frame, which is computed from the quaternion representation using the standard quaternion-to-rotation-matrix transformation, and $F$ is the total thrust generated by the four rotors. The rotational dynamics evolve according to the Euler equations of motion
\begin{align}\label{B.2}
\mathbf{I} \dot{\boldsymbol{\omega}} = \mathbf{M} - \boldsymbol{\omega} \times (\mathbf{I} \boldsymbol{\omega}),
\end{align}
with $\mathbf{I} = \text{diag}(I_{xx},I_{yy},I_{zz})$ the inertia tensor, $\boldsymbol{\omega} = [p,q,r]^{\top}$, the body angular velocity, and $\mathbf{M}=[M_{x_b},M_{y_b},M_{z_b}]^{\top}$ are the torques in the body reference frame. Expanding Eq.~\eqref{B.2} component-wise yields
\begin{align}\label{B.3}
\dot{p} &= \frac{1}{I_{xx}} \left( M_{x_b} + (I_{yy} - I_{zz})qr \right), \tag{B 3a}\\
\dot{q} &= \frac{1}{I_{yy}} \left( M_{y_b} + (I_{zz} - I_{xx})pr \right), \tag{B 3b}\\
\dot{r} &= \frac{1}{I_{zz}} \left( M_{z_b} + (I_{xx} - I_{yy})pq \right)\tag{B 3c}.
\end{align}
\addtocounter{equation}{1}

In an Euler angle representation, $\boldsymbol{\eta}=[\phi,\theta,\psi]^\top$ (roll, pitch, yaw), the angular velocities $p,q,r$ would be related to the derivatives of $\phi$, $\theta$, and $\psi$ angles through the kinematic relationship
\begin{align}\label{B.4}
\begin{bmatrix}
p \\
q \\
r
\end{bmatrix}
=
\begin{bmatrix}
\cos\theta & 0 & -\cos\phi \sin\theta \\
0 & 1 & \sin\phi \\
\sin\theta & 0 & \cos\phi \cos\theta
\end{bmatrix}
\begin{bmatrix}
\dot{\phi} \\
\dot{\theta} \\
\dot{\psi}
\end{bmatrix}.
\end{align}

However, this relationship suffers from singularities (e.g., when $\theta=\pm\pi/2$) and is not used in the implementation. Although Euler angles are commonly used for control design and output, the quadrotor attitude is represented using unit quaternions $\boldsymbol{\beta} = [\beta_0, \beta_1, \beta_2, \beta_3]^{\top}$, where $\|\boldsymbol{\beta}\| = 1$, for numerical robustness and to avoid kinematic singularities. The kinematics are given by
\begin{align}\label{B.5}
\dot{\boldsymbol{\beta}} = \frac{1}{2}
\begin{bmatrix}
0 & -p & -q & -r \\
p & 0 & r & -q \\
q & -r & 0 & p \\
r & q & -p & 0
\end{bmatrix}
\boldsymbol{\beta}.
\end{align}

For the control law, a position and attitude control structure is employed. The position loop controller computes a commanded translational acceleration by
\begin{align}\label{B.6}
\mathbf{a}_{\text{cmd}} = \mathbf{a}_{\text{des}} + \mathbf{K}_d (\dot{\mathbf{p}}_{\text{des}} - \dot{\mathbf{p}}) + \mathbf{K}_p (\mathbf{p}_{\text{des}} - \mathbf{p}),
\end{align}
where $\mathbf{K}_d$ and $\mathbf{K}_p$ are derivative and proportional gains, respectively. The commanded acceleration determines the total thrust
\begin{align}\label{B.7}
F = m(g+\mathbf{a}_{\text{cmd},z}),
\end{align}
and the desired $\phi$ and $\theta$ angles are determined as
\begin{align}\label{B.8}
\phi_{\text{des}} &= \frac{1}{g} \left( \mathbf{a}_{\text{cmd},x} \sin \psi_{\text{des}} - \mathbf{a}_{\text{cmd},y} \cos \psi_{\text{des}} \right), \tag{B 8a}\\
\theta_{\text{des}} &= \frac{1}{g} \left( \mathbf{a}_{\text{cmd},x} \cos \psi_{\text{des}} + \mathbf{a}_{\text{cmd},y} \sin \psi_{\text{des}} \right). \tag{B 8b}
\end{align}
\addtocounter{equation}{1}
The attitude control loop computes the body moments using proportional-derivative (PD) feedback by
\begin{align}\label{B.9}
\mathbf{M} = \mathbf{I}(\mathbf{K}_{dM}(\boldsymbol{\omega}_\text{des}-\boldsymbol{\omega})+\mathbf{K}_{pM}(\boldsymbol{\eta}_\text{des}-\boldsymbol{\eta})),
\end{align}
where $\boldsymbol{\omega}_{des}=[0,0,\dot\psi_{des}]^\top$, and $\mathbf{K}_{dM}$ and $\mathbf{K}_{pM}$ are attitude control gains. This controller is not used for control design in this work; its purpose is to generate informative trajectories for system identification. The resulting thrust and body moments constitute the control inputs applied to the quadrotor simulator.

To generate synthetic data of the quadrotor, we used the following physical parameters taken from \cite{LeeKimKimEtAl2025IETControlTheoryAppl}: $m = 1.3~\text{kg}$, $g = 9.81~\text{m/s}^2$, $I_{xx} = 0.0281~\text{kg}\,\text{m}^2$, $I_{yy} = 0.0286~\text{kg}\,\text{m}^2$, $I_{zz} = 0.0551~\text{kg}\,\text{m}^2$, and arm length $L=0.165~\text{m}$. Using MATLAB, we simulate the quadrotor's flight for a total duration of $30$ s using a time step of $0.01$, while tracking a reference trajectory using the PD controller \cite{MichaelMellingerLindseyEtAl2010IEEERobotAutomatMag} based on the public code of \cite{YirenLu2022}. The reference trajectory is constructed to excite the translational and rotational dynamics as
\begin{align}\label{4.11}
x(t) &= 1.5\sin (0.5t) + 0.8\sin (1.5t)+0.3\sin (3t),\tag{B 10a}\\
y(t) &= 1.5\cos (0.5t) + 0.8\cos (1.5t)+0.3\sin (3t),\tag{B 10b}\\
z(t) &= 1.5 + 0.5\sin (0.5t) + 0.3\sin (1.5t).\tag{B 10c}
\end{align}
\addtocounter{equation}{1}

\noindent with the velocities and accelerations obtained analytically by time differentiation. The $\psi(t)$ motion is independently excited via $\psi(t)=\frac{\pi}{6} \sin(t)$.

For the SI, the translational acceleration is approximated as
\begin{align}\label{B.11}
\ddot{\mathbf{x}} \approx \boldsymbol\Theta_{\text{tr}}(\mathbf{p}, \dot{\mathbf{p}}, F) \, \mathbf{W}_{\text{tr}},
\end{align}
where the library includes physics-informed terms, with polynomial terms up to second order for the velocity
\[\boldsymbol\Theta_{\text{tr}} = \left[ \mathbf{R}_{13}F \quad \mathbf{R}_{23}F \quad \mathbf{R}_{33}F \quad 1 \quad \dot{x} \quad \dot{y} \quad \dot{z} \quad \cdots \right].\]
Similarly, the rotational dynamics are identified as
\begin{equation}
\dot{\boldsymbol{\omega}} \approx \boldsymbol\Theta_{\text{ro}}(\boldsymbol{\omega}, \mathbf{M}) \, \mathbf{W}_{\text{ro}},
\end{equation}
with candidate functions
\[\boldsymbol\Theta_{\text{ro}} = \left[ M_x \quad M_y \quad M_z \quad pq \quad pr \quad qr \quad 1 \quad p \quad q \quad r \quad \cdots \right].\]

To implement MPC, we set $T_s=0.05$ s, $\Delta t^\text{M}=0.01$ s, $\Delta t^\text{plant}=0.005$ s, and $n_s^{\text{plant}}=10$. The weighting matrices $\mathbf{Q} = \text{diag}(250,250,250,1,1,1,10,10,10,10,0.1,0.1,0.1)$, prioritize position tracking (first three entries) over velocity and attitude states, with angular rates receiving the smallest penalties. The input weights $\mathbf{R}_u = \text{diag}(0.02,2,2,2)$, and $\mathbf{R}_{\Delta u} = \text{diag}(0.03, 1.5, 1.5, 1.5)$ lightly penalize thrust magnitude relative to the moments, as thrust variations are needed for altitude and obstacle avoidance maneuvers. The thrust is constrained by $F\in[0,20]$ to enforce physical limits of the motors, input $u\in[-2,1]$ reflects the physical limits of the actuators, and for the prediction and control horizons we use the values $m_p=16$ and $m_c=2$, respectively, as in \cite{LiuHongPiprekEtAl2025IEEETransAerospElectronSyst}. The safety radius $D_\text{min}=0.35 ~\text{m}$ accounts for the obstacle radius plus the drone's arm length with a small margin, and $Q_\text{obs}=1500$. At each prediction step, the reference state includes the desired position trajectory, the quaternion reference is fixed to the identity quaternion $[1,0,0,0]$ to encourage upright flight, and the observed quaternion is normalized to unit norm.

\begin{figure}[tb]
\centering\includegraphics[width=5.0in]{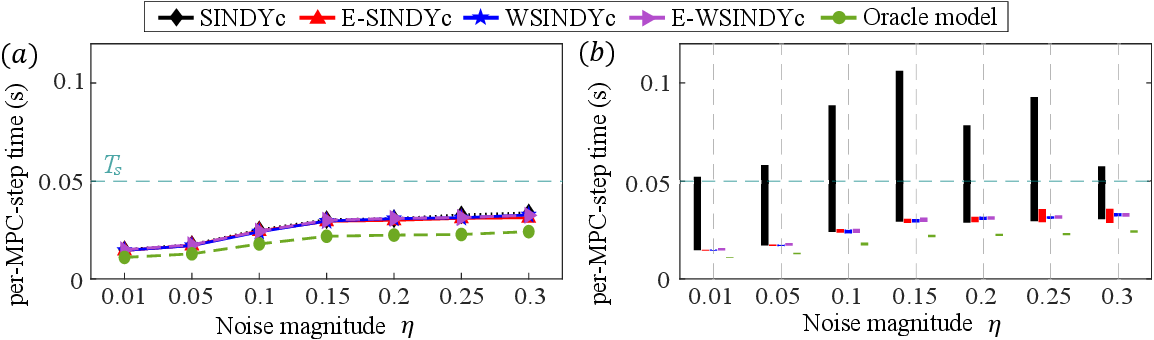}
\caption{Drone, per-MPC-step runtime using Eq.~\eqref{5.1}: ($a$) Medians ($b$) first to third quartile of ($a$), from 200 noise realizations each (Online version in color.)}
\label{fig_drone_MPC_time}
\vspace{-7mm}
\end{figure}

The per-MPC-step solve time is shown in Fig.~\ref{fig_drone_MPC_time}. Except for the SINDY--MPC framework, the other sparse-regression-based methods ensure that in all cases the computation terminates on schedule.

\begin{figure}[tb]
\centering\includegraphics[width=5.0in]{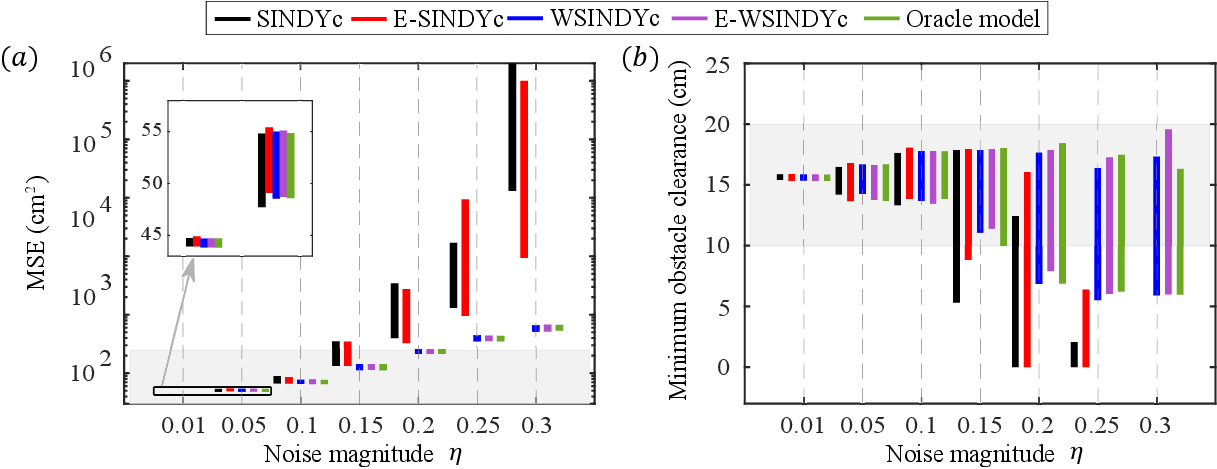}
\caption{Trajectory tracking performance and obstacle avoidance. Performance for different noise levels: ($a$) mean square error, and ($b$) minimum obstacle clearance, for 200 noise realizations each. Outside the shaded region, prediction accuracy deteriorates, and the resulting control actions fail to guarantee safe obstacle avoidance.  (Online version in color.)}
\label{fig_drone_dep_box}
\vspace{-7mm}
\end{figure}

\section*{Appendix C. Lorenz system and simulation parameters}
\refstepcounter{section}\label{appendix:lorenz}
The governing equations of this system are given by:
\begin{align}\
\dot{x}_1 &= 10 (x_2 - x_1) + u, \tag{C 1a}\label{C1a} \\
\dot{x}_2 &= x_1(28-x_3) - x_2, \tag{C 1b}\label{C1b} \\
\dot{x}_3 &= x_1x_2 - \frac{8}{3} x_3. \tag{C 1c}\label{C1c}
\end{align}
\addtocounter{equation}{1}

To generate the synthetic time series, MATLAB's \texttt{ode45} solver is used, with the relative and absolute tolerances set to $10^{-10}$. For this system, the Schroeder sweep signal \cite{Schroeder1970IEEETransInformTheory} is applied as the actuation input $u(t)$. The DMDc is trained on data shifted by the reference state, so that the model describes the evolution of deviations from the goal, and the reference state is added back after predictions. For the NN, $1$ hidden layer with $10$ neurons is used. 
In the MPC implementation, we set $T_s=0.01$ s, $\Delta t^\text{M}=\Delta t^\text{plant}=0.001$ s, and $n_s^{\text{plant}}=10$. The target point is $(-\sqrt{72},-\sqrt{72},27)$. The weighting matrices are $\mathbf{Q} = \text{diag}(1,1,1)$ and $\mathbf{R}_u = \mathbf{R}_{\Delta u} = 0.001$, the control input is constrained to $u\in[-50,50]$, and $m_p=m_c=10$.

\begin{figure}[tb]
\centering\includegraphics[width=5in]{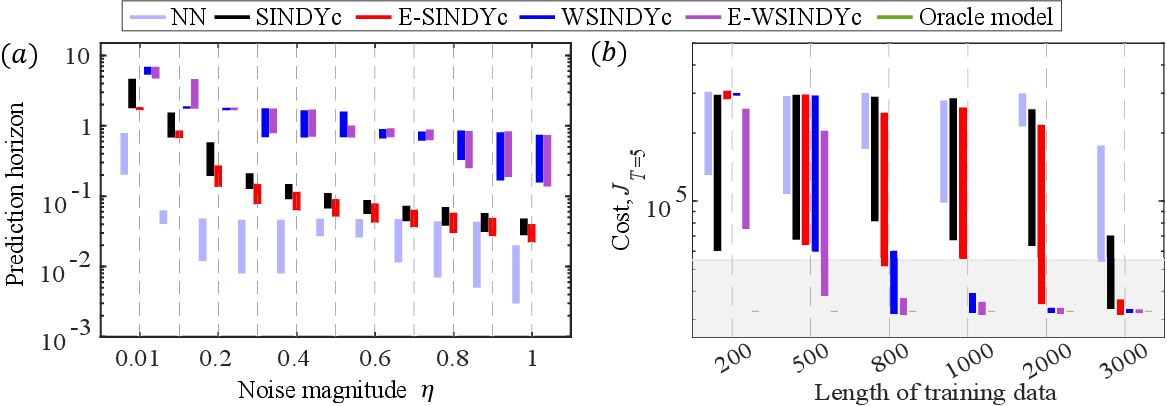}
\caption{Lorenz System, first to third quartile: ($a$) Validation prediction performance for different noise levels, prediction horizon in time units, ($b$) different training data, terminal cumulative cost, using Eq.~\eqref{3.1}, over 5 time units, from 200 noise realizations each (Online version in color.)}
\label{fig_Lorenz2_box}
\vspace{-7mm}
\end{figure}

\begin{figure}[tb]
\centering\includegraphics[width=5in]{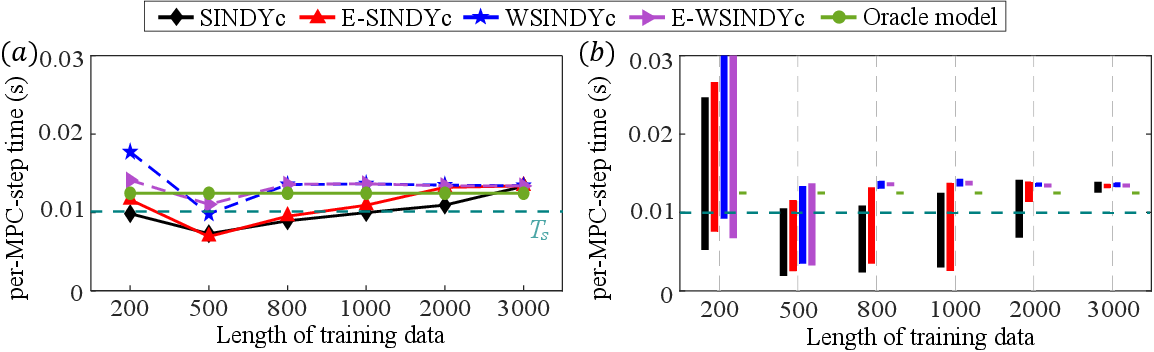}
\caption{Lorenz system, per-MPC-step optimization time Eq.~\eqref{3.1}: ($a$) Medians ($b$) first to third quartile of ($a$), from 200 noise realizations each (Online version in color.)}
\label{fig_Lorenz4}
\vspace{-7mm}
\end{figure}

Figure ~\ref{fig_Lorenz4} shows the per-MPC-step solve times corresponding to Fig.~\ref{fig_Lorenz3}d, obtained using the C implementation. As in the previous two systems, this implementation accelerates the approach of \cite{KaiserKutzBrunton2018ProcRSocA}, bringing the solve times close to $T_s$.

\section*{Appendix D. MPC on a F-8 aircraft}
\refstepcounter{section}\label{appendix:f8}
This section considers an automatic flight control system for the F-8 Crusader operating at 9000 m (30000 ft) and a Mach number of 0.85\cite{GarrardJordan1977Automatica,YanWang2012IEEETransIndInf}, where the objective is to track a specified angle of attack trajectory. The aircraft dynamics are given by:
\begin{align}
\dot{x}_1 &= 
\begin{aligned}[t]
&-0.877x_1 + x_3 - 0.088x_1x_3 + 0.47x_1^2 - 0.019x_2^2 - x_1^2x_3 + 3.846x_1^3 - 0.215u \\
&+ 0.28x_1^2u + 0.47x_1u^2 + 0.63u^3,
\end{aligned}
\tag{D 1a}\label{D1a} \\
\dot{x}_2 &= x_3,
\tag{D 1b}\label{D1b} \\
\dot{x}_3 &= 
\begin{aligned}[t]
-4.208x_1 - 0.396x_3 - 0.47x_1^2 - 3.564x_1^3 - 20.967u + 6.265x_1^2u + 46x_1u^2 + 61.4u^3,
\end{aligned}
\tag{D 1c}\label{D1c} 
\end{align}\addtocounter{equation}{1}
where $x_1$, $x_2$, $x_3$, and $u$ are the angle of attack, pitch angle, pitch rate, and tail deflection angle. The system is highly nonlinear due to non-affine coupling between states and the control input. The angle of attack to be tracked is:
\begin{align}
r=0.4 \left( 
        -\frac{0.5}{1 + \text{e}^{\frac{t}{0.1} - 0.8}} 
        + \frac{1}{1 + \text{e}^{\frac{t}{0.1} - 3}} 
        - 0.4 \right).
\tag{D 2}\label{D2}
\end{align} 
The control objective focuses exclusively on tracking the angle of attack $x_1$, which is designated as the output variable $y$ for the performance optimization. The control inputs are recalculated every 10 time steps and held constant in between. For MPC, $T_s=0.01$ s, the model timestep is $\Delta t^\text{M}=0.01$ s and the plant's $\Delta t^\text{plant}=0.001$ s, and $n_s^{\text{plant}}=10$. The weighting matrices are $\mathbf{Q} = \text{diag}(25,0,0)$ and $\mathbf{R}_u = \mathbf{R}_{\Delta u} = 0.05$, with input constrains to $u\in[-0.3,0.5]$, $\Delta u \in[-0.1,0.1]$, output limits $y\in[-0.2,0.4]$, and $m_p=m_c=13$.

\begin{figure}[tb]
\centering\includegraphics[width=5in]{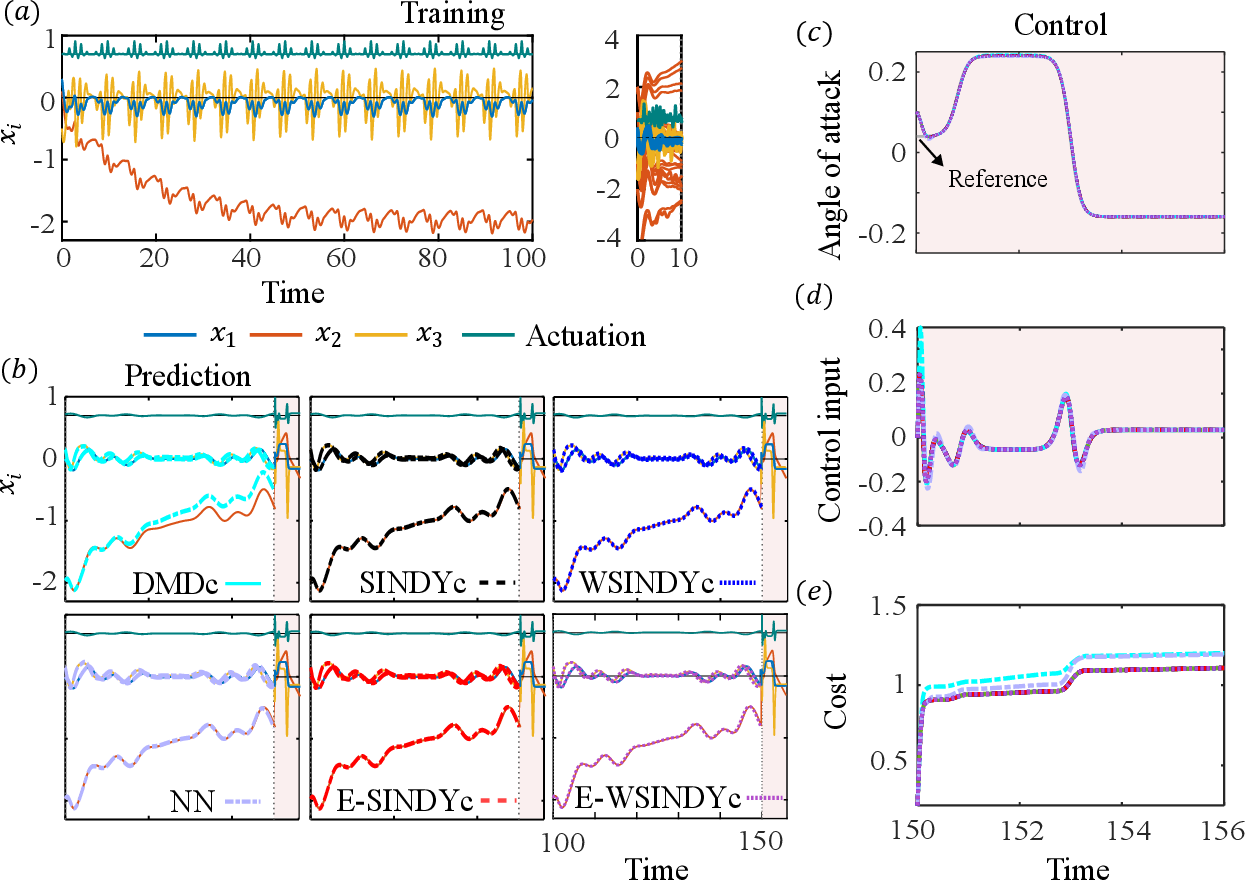}
\caption{Prediction and control performance for the F-8 aircraft model: ($a$) time series of states and input during training: left panel for DMDc and SINDYc-based methods; right panel for NN, ($b$)  Validation stage, ($c$) Reference ($r$) and angle of attack $y=x_1$, ($d$) control input, and ($e$) cumulative cost over 6 time units during control. (Online version in color.)}
\label{fig_F8}
\vspace{-7mm}
\end{figure}

In this case, the NN has 2 hidden layers, each with 15 neurons. To obtain desirable results, it was trained using 250 different short trajectories of 1000 snapshots each, of which 20 samples are plotted in the right panel of Fig.~\ref{fig_F8}a. From Fig.~\ref{fig_F8}b, except for DMDc, the other methods accurately discover the relationships between variables. All methods perform well for MPC (Fig.~\ref{fig_F8}c-e), but the NN achieves this using a much larger amount of data. 

\begin{figure}[tb]
\centering\includegraphics[width=3.5in]{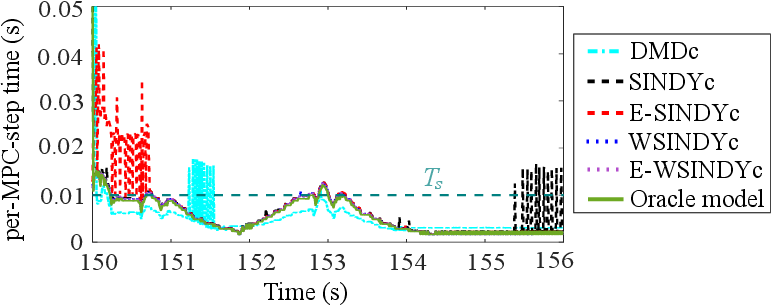}
\caption{F-8 aircraft model, per-MPC-step runtime using Eq.~\eqref{3.1} (Online version in color.)}
\label{fig_F8_1}
\vspace{-7mm}
\end{figure}

In Fig. ~\ref{fig_F8_1}, we display the per-MPC-step solve times associated to Fig.~\ref{fig_F8}c-d, employing the C implementation. As can be seen, for all models, the per-MPC-step times are mostly below $T_s$.

\vskip2pc

\bibliography{WSINDyControl}

@article{AbdullahChristofides2023ComputChemEng,
  title = {Data-Based Modeling and Control of Nonlinear Process Systems Using Sparse Identification: {{An}} Overview of Recent Results},
  shorttitle = {Data-Based Modeling and Control of Nonlinear Process Systems Using Sparse Identification},
  author = {Abdullah, Fahim and Christofides, Panagiotis D.},
  year = 2023,
  month = jun,
  journal = {Computers \& Chemical Engineering},
  volume = {174},
  pages = {108247},
  issn = {00981354},
  doi = {10.1016/j.compchemeng.2023.108247},
  urldate = {2023-07-13},
  langid = {english}
}

@article{AbdullahChristofides2023ChemicalEngineeringResearchandDesign,
  title = {Real-Time Adaptive Sparse-Identification-Based Predictive Control of Nonlinear Processes},
  author = {Abdullah, Fahim and Christofides, Panagiotis D.},
  year = 2023,
  month = aug,
  journal = {Chemical Engineering Research and Design},
  volume = {196},
  pages = {750--769},
  issn = {02638762},
  doi = {10.1016/j.cherd.2023.07.011},
  urldate = {2026-06-01},
  langid = {english}
}

@article{AbdullahAlhajeriChristofides2022IndEngChemRes,
  title = {Modeling and {{Control}} of {{Nonlinear Processes Using Sparse Identification}}: {{Using Dropout}} to {{Handle Noisy Data}}},
  shorttitle = {Modeling and {{Control}} of {{Nonlinear Processes Using Sparse Identification}}},
  author = {Abdullah, Fahim and Alhajeri, Mohammed S. and Christofides, Panagiotis D.},
  year = 2022,
  month = dec,
  journal = {Industrial \& Engineering Chemistry Research},
  volume = {61},
  number = {49},
  pages = {17976--17992},
  issn = {0888-5885, 1520-5045},
  doi = {10.1021/acs.iecr.2c02639},
  urldate = {2026-06-01},
  copyright = {https://doi.org/10.15223/policy-029},
  langid = {english}
}

@article{AbdullahWuChristofides2021ComputChemEng,
  ids = {AbdullahWuChristofides2021ComputersChemicalEngineering},
  title = {Sparse-Identification-Based Model Predictive Control of Nonlinear Two-Time-Scale Processes},
  author = {Abdullah, Fahim and Wu, Zhe and Christofides, Panagiotis D.},
  year = 2021,
  month = jun,
  journal = {Computers \& Chemical Engineering},
  pages = {107411},
  issn = {00981354},
  doi = {10.1016/j.compchemeng.2021.107411},
  urldate = {2021-06-23},
  langid = {english}
}

@article{Akaike1974IEEETransAutomControl,
  title = {A New Look at the Statistical Model Identification},
  author = {Akaike, H.},
  year = 1974,
  month = dec,
  journal = {IEEE Transactions on Automatic Control},
  volume = {19},
  number = {6},
  pages = {716--723},
  issn = {0018-9286},
  doi = {10.1109/TAC.1974.1100705},
  urldate = {2014-11-12},
  langid = {english}
}

@article{BerkoozHolmesLumley1993AnnuRevFluidMech,
  title = {The {{Proper Orthogonal Decomposition}} in the {{Analysis}} of {{Turbulent Flows}}},
  author = {Berkooz, G and Holmes, P and Lumley, J L},
  year = 1993,
  month = jan,
  journal = {Annual Review of Fluid Mechanics},
  volume = {25},
  number = {1},
  pages = {539--575},
  issn = {0066-4189, 1545-4479},
  doi = {10.1146/annurev.fl.25.010193.002543},
  urldate = {2023-10-03},
  langid = {english}
}

@article{BiekerPeitzBruntonEtAl2020TheorComputFluidDyn,
  title = {Deep Model Predictive Flow Control with Limited Sensor Data and Online Learning},
  author = {Bieker, Katharina and Peitz, Sebastian and Brunton, Steven L. and Kutz, J. Nathan and Dellnitz, Michael},
  year = 2020,
  month = aug,
  journal = {Theoretical and Computational Fluid Dynamics},
  volume = {34},
  number = {4},
  pages = {577--591},
  issn = {0935-4964, 1432-2250},
  doi = {10.1007/s00162-020-00520-4},
  urldate = {2025-10-03},
  langid = {english}
}

@article{BoninsegnaNuskeClementi2018JChemPhys,
  ids = {BoninsegnaNuskeClementi2018TheJournalofChemicalPhysicsa},
  title = {Sparse Learning of Stochastic Dynamical Equations},
  author = {Boninsegna, Lorenzo and N{\"u}ske, Feliks and Clementi, Cecilia},
  year = 2018,
  month = jun,
  journal = {The Journal of Chemical Physics},
  volume = {148},
  number = {24},
  pages = {241723},
  issn = {0021-9606, 1089-7690},
  doi = {10.1063/1.5018409},
  urldate = {2020-06-11},
  langid = {english}
}

@article{BonninDekeyserPittsEtAl2016PlasmaandFusionResearch,
  title = {Presentation of the {{New SOLPS-ITER Code Package}} for {{Tokamak Plasma Edge Modelling}}},
  author = {Bonnin, Xavier and Dekeyser, Wouter and Pitts, Richard and Coster, David and Voskoboynikov, Serguey and Wiesen, Sven},
  year = 2016,
  journal = {Plasma and Fusion Research},
  volume = {11},
  number = {0},
  pages = {1403102--1403102},
  issn = {1880-6821},
  doi = {10.1585/pfr.11.1403102},
  urldate = {2026-04-20},
  langid = {english}
}

@article{BoraseMaghadeSondkarEtAl2021IntJDynamControl,
  title = {A Review of {{PID}} Control, Tuning Methods and Applications},
  author = {Borase, Rakesh P. and Maghade, D. K. and Sondkar, S. Y. and Pawar, S. N.},
  year = 2021,
  month = jun,
  journal = {International Journal of Dynamics and Control},
  volume = {9},
  number = {2},
  pages = {818--827},
  issn = {2195-268X, 2195-2698},
  doi = {10.1007/s40435-020-00665-4},
  urldate = {2026-02-13},
  langid = {english}
}

@article{BortzMessengerDukic2023BullMathBiol,
  title = {Direct {{Estimation}} of {{Parameters}} in {{ODE Models Using WENDy}}: {{Weak-form Estimation}} of {{Nonlinear Dynamics}}},
  author = {Bortz, David M. and Messenger, Daniel A. and Dukic, Vanja},
  year = 2023,
  journal = {Bulletin of Mathematical Biology},
  volume = {85},
  number = {110},
  doi = {10.1007/S11538-023-01208-6}
}

@incollection{BortzMessengerTran2024NumericalAnalysisMeetsMachineLearning,
  title = {Weak Form-Based Data-Driven Modeling: {{Computationally Efficient}} and {{Noise Robust Equation Learning}} and {{Parameter Inference}}},
  booktitle = {Numerical {{Analysis Meets Machine Learning}}},
  author = {Bortz, David M. and Messenger, Daniel A. and Tran, April},
  editor = {Mishra, Siddhartha and Townsend, Alex},
  year = 2024,
  series = {Handbook of {{Numerical Analysis}}},
  volume = {25},
  pages = {54--82},
  publisher = {Elsevier}
}

@article{Bottcher2026NonlinearDyn,
  title = {Control of Dynamical Systems with Neural Networks},
  author = {B{\"o}ttcher, Lucas},
  year = 2026,
  month = jan,
  journal = {Nonlinear Dynamics},
  volume = {114},
  number = {2},
  pages = {79},
  issn = {0924-090X, 1573-269X},
  doi = {10.1007/s11071-025-11937-z},
  urldate = {2026-02-13},
  langid = {english}
}

@book{BruntonKutz2019,
  title = {Data-{{Driven Science}} and {{Engineering}}: {{Machine Learning}}, {{Dynamical Systems}}, and {{Control}}},
  shorttitle = {Data-{{Driven Science}} and {{Engineering}}},
  author = {Brunton, Steven L. and Kutz, J. Nathan},
  year = 2019,
  month = jan,
  edition = {1},
  publisher = {Cambridge University Press},
  doi = {10.1017/9781108380690},
  urldate = {2019-06-07},
  isbn = {978-1-108-38069-0 978-1-108-42209-3}
}

@article{BruntonProctorKutz2016IFAC-PapersOnLine,
  title = {Sparse {{Identification}} of {{Nonlinear Dynamics}} with {{Control}} ({{SINDYc}})},
  author = {Brunton, Steven L. and Proctor, Joshua L. and Kutz, J. Nathan},
  year = 2016,
  journal = {IFAC-PapersOnLine},
  volume = {49},
  number = {18},
  pages = {710--715},
  issn = {24058963},
  doi = {10.1016/j.ifacol.2016.10.249},
  urldate = {2021-04-16},
  langid = {english}
}

@article{BruntonProctorKutz2016ProcNatlAcadSci,
  title = {Discovering Governing Equations from Data by Sparse Identification of Nonlinear Dynamical Systems},
  author = {Brunton, Steven L. and Proctor, Joshua L. and Kutz, J. Nathan},
  year = 2016,
  month = apr,
  journal = {Proceedings of the National Academy of Sciences},
  volume = {113},
  number = {15},
  pages = {3932--3937},
  issn = {0027-8424, 1091-6490},
  doi = {10.1073/pnas.1517384113},
  urldate = {2019-12-31},
  langid = {english}
}

@article{BruntonZolmanKutzEtAl2025AnnuRevControlRobotAutonSyst,
  title = {Machine {{Learning}} for {{Sparse Nonlinear Modeling}} and {{Control}}},
  author = {Brunton, Steven L. and Zolman, Nicholas and Kutz, J. Nathan and Fasel, Urban},
  year = 2025,
  month = jan,
  journal = {Annual Review of Control, Robotics, and Autonomous Systems},
  issn = {2573-5144},
  doi = {10.1146/annurev-control-030123-015238},
  urldate = {2025-03-13},
  langid = {english}
}

@article{Camps-VallsGerhardusNinadEtAl2023PhysicsReports,
  title = {Discovering Causal Relations and Equations from Data},
  author = {{Camps-Valls}, Gustau and Gerhardus, Andreas and Ninad, Urmi and Varando, Gherardo and Martius, Georg and {Balaguer-Ballester}, Emili and Vinuesa, Ricardo and Diaz, Emiliano and Zanna, Laure and Runge, Jakob},
  year = 2023,
  month = dec,
  journal = {Physics Reports},
  volume = {1044},
  pages = {1--68},
  issn = {03701573},
  doi = {10.1016/j.physrep.2023.10.005},
  urldate = {2023-11-19},
  langid = {english}
}

@article{CaoLaiAlam2017JIntellRobotSyst,
  title = {Gaussian {{Process Model Predictive Control}} of an {{Unmanned Quadrotor}}},
  author = {Cao, Gang and Lai, Edmund M.-K. and Alam, Fakhrul},
  year = 2017,
  month = oct,
  journal = {Journal of Intelligent \& Robotic Systems},
  volume = {88},
  number = {1},
  pages = {147--162},
  issn = {0921-0296, 1573-0409},
  doi = {10.1007/s10846-017-0549-y},
  urldate = {2026-06-13},
  langid = {english}
}

@article{CenedeseAxasBauerleinEtAl2022NatCommun,
  title = {Data-Driven Modeling and Prediction of Non-Linearizable Dynamics via Spectral Submanifolds},
  author = {Cenedese, Mattia and Ax{\aa}s, Joar and B{\"a}uerlein, Bastian and Avila, Kerstin and Haller, George},
  year = 2022,
  month = feb,
  journal = {Nature Communications},
  volume = {13},
  number = {1},
  pages = {872},
  issn = {2041-1723},
  doi = {10.1038/s41467-022-28518-y},
  urldate = {2023-03-26},
  langid = {english}
}

@misc{ChenWang2025arXiv250906882,
  title = {Dynamic {{Modeling}} and {{Efficient Data-Driven Optimal Control}} for {{Micro Autonomous Surface Vehicles}}},
  author = {Chen, Zhiheng and Wang, Wei},
  year = 2025,
  month = sep,
  number = {arXiv:2509.06882},
  eprint = {2509.06882},
  primaryclass = {cs},
  publisher = {arXiv},
  doi = {10.48550/arXiv.2509.06882},
  urldate = {2025-09-09},
  archiveprefix = {arXiv},
  langid = {english}
}

@inproceedings{ChenRubanovaBettencourtEtAl2018AdvancesinNeuralInformationProcessingSystems,
  title     = {Neural Ordinary Differential Equations},
  author    = {Chen, Tian Qi and Rubanova, Yulia and Bettencourt, Jesse and Duvenaud, David K},
  booktitle = {Advances in Neural Information Processing Systems},
  year      = {2018},
  volume    = {31}
}

@article{CrutchfieldMcNamara1987ComplexSyst,
  title = {Equations of {{Motion}} from a {{Data Series}}},
  author = {Crutchfield, James P and McNamara, Bruce S.},
  year = 1987,
  journal = {Complex Systems},
  volume = {1},
  number = {3},
  pages = {417--452},
  langid = {english}
}

@inproceedings{FaselKaiserKutzEtAl2021202160thIEEEConfDecisControlCDC,
  title = {{{SINDy}} with {{Control}}: {{A Tutorial}}},
  shorttitle = {{{SINDy}} with {{Control}}},
  booktitle = {2021 60th {{IEEE Conference}} on {{Decision}} and {{Control}} ({{CDC}})},
  author = {Fasel, Urban and Kaiser, Eurika and Kutz, J. Nathan and Brunton, Bingni W. and Brunton, Steven L.},
  year = 2021,
  month = dec,
  pages = {16--21},
  publisher = {IEEE},
  address = {Austin, TX, USA},
  doi = {10.1109/CDC45484.2021.9683120},
  urldate = {2022-05-03},
  isbn = {978-1-6654-3659-5},
  langid = {english}
}

@article{FaselKutzBruntonEtAl2022ProcRSocMathPhysEngSci,
  title = {Ensemble-{{SINDy}}: {{Robust}} Sparse Model Discovery in the Low-Data, High-Noise Limit, with Active Learning and Control},
  shorttitle = {Ensemble-{{SINDy}}},
  author = {Fasel, U. and Kutz, J. N. and Brunton, B. W. and Brunton, S. L.},
  year = 2022,
  month = apr,
  journal = {Proceedings of the Royal Society A: Mathematical, Physical and Engineering Sciences},
  volume = {478},
  number = {2260},
  pages = {20210904},
  issn = {1364-5021, 1471-2946},
  doi = {10.1098/rspa.2021.0904},
  urldate = {2022-04-20},
  langid = {english}
}

@article{FreyBaumgartnerFrisonEtAl2025OptimControlApplMethods,
  title = {Multi-{{Phase Optimal Control Problems}} for {{Efficient Nonlinear Model Predictive Control With}} Acados},
  author = {Frey, Jonathan and Baumg{\"a}rtner, Katrin and Frison, Gianluca and Diehl, Moritz},
  year = 2025,
  month = mar,
  journal = {Optimal Control Applications and Methods},
  volume = {46},
  number = {2},
  pages = {827--845},
  issn = {0143-2087, 1099-1514},
  doi = {10.1002/oca.3234},
  urldate = {2026-06-10},
  langid = {english}
}

@article{FukamiMurataZhangEtAl2021JFluidMech,
  title = {Sparse Identification of Nonlinear Dynamics with Low-Dimensionalized Flow Representations},
  author = {Fukami, Kai and Murata, Takaaki and Zhang, Kai and Fukagata, Koji},
  year = 2021,
  month = nov,
  journal = {Journal of Fluid Mechanics},
  volume = {926},
  pages = {A10},
  issn = {0022-1120, 1469-7645},
  doi = {10.1017/jfm.2021.697},
  urldate = {2021-11-10},
  langid = {english}
}

@article{GarciaPrettMorari1989Automatica,
  title = {Model Predictive Control: {{Theory}} and Practice---{{A}} Survey},
  shorttitle = {Model Predictive Control},
  author = {Garc{\'i}a, Carlos E. and Prett, David M. and Morari, Manfred},
  year = 1989,
  month = may,
  journal = {Automatica},
  volume = {25},
  number = {3},
  pages = {335--348},
  issn = {00051098},
  doi = {10.1016/0005-1098(89)90002-2},
  urldate = {2025-09-30},
  copyright = {https://www.elsevier.com/tdm/userlicense/1.0/},
  langid = {english}
}

@article{GarrardJordan1977Automatica,
  title = {Design of Nonlinear Automatic Flight Control Systems},
  author = {Garrard, William L. and Jordan, John M.},
  year = 1977,
  month = sep,
  journal = {Automatica},
  volume = {13},
  number = {5},
  pages = {497--505},
  issn = {00051098},
  doi = {10.1016/0005-1098(77)90070-X},
  urldate = {2025-10-10},
  copyright = {https://www.elsevier.com/tdm/userlicense/1.0/},
  langid = {english}
}

@article{GuevaraVarela-AldasGandolfoEtAl2025Drones,
  title = {{{SINDy}} and {{PD-Based UAV Dynamics Identification}} for {{MPC}}},
  author = {Guevara, Bryan S. and {Varela-Ald{\'a}s}, Jos{\'e} and Gandolfo, Daniel C. and Toibero, Juan M.},
  year = 2025,
  month = jan,
  journal = {Drones},
  volume = {9},
  number = {1},
  pages = {71},
  issn = {2504-446X},
  doi = {10.3390/drones9010071},
  urldate = {2026-01-08},
  langid = {english}
}

@article{GurevichReinboldGrigoriev2019Chaos,
  title = {Robust and Optimal Sparse Regression for Nonlinear {{PDE}} Models},
  author = {Gurevich, Daniel R. and Reinbold, Patrick A. K. and Grigoriev, Roman O.},
  year = 2019,
  month = oct,
  journal = {Chaos: An Interdisciplinary Journal of Nonlinear Science},
  volume = {29},
  number = {10},
  pages = {103113},
  issn = {1054-1500, 1089-7682},
  doi = {10.1063/1.5120861},
  urldate = {2021-01-05},
  langid = {english}
}

@inproceedings{HaoYaoSuEtAl2024AdvNeuralInfProcessSyst,
  title = {{{PINNacle}}: {{A Comprehensive Benchmark}} of {{Physics-Informed Neural Networks}} for {{Solving PDEs}}},
  booktitle = {Advances in {{Neural Information Processing Systems}}},
  author = {Hao, Zhongkai and Yao, Jiachen and Su, Chang and Su, Hang and Wang, Ziao and Lu, Fanzhi and Xia, Zeyu and Zhang, Yichi and Liu, Songming and Lu, Lu and Zhu, Jun},
  editor = {Globerson, A. and Mackey, L. and Belgrave, D. and Fan, A. and Paquet, U. and Tomczak, J. and Zhang, C.},
  year = 2024,
  volume = {37},
  pages = {76721--76774},
  publisher = {Curran Associates, Inc.},
  langid = {english}
}

@article{HeydariNovinzadehTayefi2025AdaptiveControlampSignal,
  title = {Anti {{Wind}}-{{Up}} and {{Robust Data}}-{{Driven Model}}-{{Free Adaptive Control}} for {{MIMO Nonlinear Discrete}}-{{Time Systems}}},
  author = {Heydari, Mohsen and Novinzadeh, Alireza B. and Tayefi, Morteza},
  year = 2025,
  month = oct,
  journal = {International Journal of Adaptive Control and Signal Processing},
  volume = {39},
  number = {10},
  pages = {2143--2157},
  issn = {0890-6327, 1099-1115},
  doi = {10.1002/acs.3907},
  urldate = {2026-06-10},
  langid = {english}
}

@article{JiangXiongZhangEtAl2021NonlinearDyn,
  title = {Modeling and Prediction of the Transmission Dynamics of {{COVID-19}} Based on the {{SINDy-LM}} Method},
  author = {Jiang, Yu-Xin and Xiong, Xiong and Zhang, Shuo and Wang, Jia-Xiang and Li, Jia-Chun and Du, Lin},
  year = 2021,
  month = aug,
  journal = {Nonlinear Dynamics},
  volume = {105},
  number = {3},
  pages = {2775--2794},
  issn = {0924-090X, 1573-269X},
  doi = {10.1007/s11071-021-06707-6},
  urldate = {2022-01-23},
  langid = {english}
}

@article{JuangPappa1985JournalofGuidanceControlandDynamics,
  title = {An Eigensystem Realization Algorithm for Modal Parameter Identification and Model Reduction},
  author = {Juang, Jer-Nan and Pappa, Richard S.},
  year = 1985,
  month = sep,
  journal = {Journal of Guidance, Control, and Dynamics},
  volume = {8},
  number = {5},
  pages = {620--627},
  issn = {0731-5090, 1533-3884},
  doi = {10.2514/3.20031},
  urldate = {2025-09-29},
  langid = {english}
}

@article{KaiserKutzBrunton2018ProcRSocA,
  ids = {KaiserKutzBrunton2018ProcRSocAa},
  title = {Sparse Identification of Nonlinear Dynamics for Model Predictive Control in the Low-Data Limit},
  author = {Kaiser, E. and Kutz, J. N. and Brunton, S. L.},
  year = 2018,
  month = nov,
  journal = {Proceedings of the Royal Society A: Mathematical, Physical and Engineering Sciences},
  volume = {474},
  number = {2219},
  pages = {20180335},
  issn = {1364-5021, 1471-2946},
  doi = {10.1098/rspa.2018.0335},
  urldate = {2021-04-16},
  langid = {english}
}

@article{KaptanogluHansenLoreEtAl2023PhysPlasmas,
  title = {Sparse Regression for Plasma Physics},
  author = {Kaptanoglu, Alan A. and Hansen, Christopher and Lore, Jeremy D. and Landreman, Matt and Brunton, Steven L.},
  year = 2023,
  month = mar,
  journal = {Physics of Plasmas},
  volume = {30},
  number = {3},
  pages = {033906},
  issn = {1070-664X, 1089-7674},
  doi = {10.1063/5.0139039},
  urldate = {2025-08-20},
  langid = {english}
}

@article{KordaMezic2018Automatica,
  title = {Linear Predictors for Nonlinear Dynamical Systems: {{Koopman}} Operator Meets Model Predictive Control},
  shorttitle = {Linear Predictors for Nonlinear Dynamical Systems},
  author = {Korda, Milan and Mezi{\'c}, Igor},
  year = 2018,
  month = jul,
  journal = {Automatica},
  volume = {93},
  pages = {149--160},
  issn = {00051098},
  doi = {10.1016/j.automatica.2018.03.046},
  urldate = {2026-02-20},
  langid = {english}
}

@article{LaiNagarajaiah2019MechSystSignalProcess,
  title = {Sparse Structural System Identification Method for Nonlinear Dynamic Systems with Hysteresis/Inelastic Behavior},
  author = {Lai, Zhilu and Nagarajaiah, Satish},
  year = 2019,
  month = feb,
  journal = {Mechanical Systems and Signal Processing},
  volume = {117},
  pages = {813--842},
  issn = {08883270},
  doi = {10.1016/j.ymssp.2018.08.033},
  urldate = {2022-02-05},
  langid = {english}
}

@article{LeeCarlberg2020JComputPhys,
  title = {Model Reduction of Dynamical Systems on Nonlinear Manifolds Using Deep Convolutional Autoencoders},
  author = {Lee, Kookjin and Carlberg, Kevin T.},
  year = 2020,
  month = mar,
  journal = {Journal of Computational Physics},
  volume = {404},
  pages = {108973},
  issn = {00219991},
  doi = {10.1016/j.jcp.2019.108973},
  urldate = {2021-04-25},
  langid = {english}
}

@article{LeeKimKimEtAl2025IETControlTheoryAppl,
  title = {Sparse {{Identification}} of {{Nonlinear Dynamics}}-{{Based Model Predictive Control}} for {{Multirotor Collision Avoidance}}},
  author = {Lee, Jayden Dongwoo and Kim, Youngjae and Kim, Yoonseong and Bang, Hyochoong},
  year = 2025,
  month = jan,
  journal = {IET Control Theory \& Applications},
  volume = {19},
  number = {1},
  pages = {e70049},
  issn = {1751-8644, 1751-8652},
  doi = {10.1049/cth2.70049},
  urldate = {2026-01-08},
  langid = {english}
}

@article{LeeRenQianEtAl2025IEEEASMETransMechatron,
  title = {Energy {{Reduction}} for {{Wearable Pneumatic Valve System With SINDy}} and {{Time-Variant Model Predictive Control}}},
  author = {Lee, Hao and Ren, Ruoning and Qian, Yifei and Rosen, Jacob},
  year = 2025,
  month = apr,
  journal = {IEEE/ASME Transactions on Mechatronics},
  volume = {30},
  number = {2},
  pages = {862--872},
  publisher = {{Institute of Electrical and Electronics Engineers (IEEE)}},
  issn = {1083-4435, 1941-014X},
  doi = {10.1109/tmech.2024.3458092},
  urldate = {2025-07-27},
  copyright = {https://ieeexplore.ieee.org/Xplorehelp/downloads/license-information/IEEE.html},
  langid = {english}
}

@article{LejarzaBaldea2022SciRep,
  title = {Data-Driven Discovery of the Governing Equations of Dynamical Systems via Moving Horizon Optimization},
  author = {Lejarza, Fernando and Baldea, Michael},
  year = 2022,
  month = jul,
  journal = {Scientific Reports},
  volume = {12},
  number = {1},
  pages = {11836},
  issn = {2045-2322},
  doi = {10.1038/s41598-022-13644-w},
  urldate = {2023-09-16},
  langid = {english}
}

@article{LiuHongPiprekEtAl2025IEEETransAerospElectronSyst,
  title = {Data-{{Based Modeling}} and {{Control}} of the {{Nonlinear Aircraft System Using Extended Implicit Sparse Identification}}},
  author = {Liu, Yue and Hong, Haichao and Piprek, Patrick and Chud{\'y}, Peter and Hu, Shiqiang},
  year = 2025,
  month = jun,
  journal = {IEEE Transactions on Aerospace and Electronic Systems},
  volume = {61},
  number = {3},
  pages = {6928--6940},
  issn = {0018-9251, 1557-9603, 2371-9877},
  doi = {10.1109/TAES.2025.3531345},
  urldate = {2026-01-08},
  copyright = {https://ieeexplore.ieee.org/Xplorehelp/downloads/license-information/IEEE.html},
  langid = {english}
}

@article{LiuMaWangEtAl2024arXiv240810205,
  title = {{{KAN}} 2.0: {{Kolmogorov-Arnold Networks Meet Science}}},
  shorttitle = {{{KAN}} 2.0},
  author = {Liu, Ziming and Ma, Pingchuan and Wang, Yixuan and Matusik, Wojciech and Tegmark, Max},
  year = 2024,
  month = aug,
  eprint = {2408.10205},
  primaryclass = {physics},
  publisher = {arXiv},
  urldate = {2024-09-21},
  archiveprefix = {arXiv},
  langid = {english}
}

@book{Ljung1999,
  title = {System Identification: Theory for the User},
  shorttitle = {System Identification},
  author = {Ljung, Lennart},
  year = 1999,
  edition = {2nd},
  publisher = {Prentice Hall PTR},
  address = {Upper Saddle River, NJ},
  isbn = {978-0-13-656695-3},
  lccn = {QA402 .L59 1999}
}

@article{Lopez-SanchezMoreno-Valenzuela2023AnnualReviewsinControl,
  title = {{{PID}} Control of Quadrotor {{UAVs}}: {{A}} Survey},
  shorttitle = {{{PID}} Control of Quadrotor {{UAVs}}},
  author = {{Lopez-Sanchez}, Ivan and {Moreno-Valenzuela}, Javier},
  year = 2023,
  journal = {Annual Reviews in Control},
  volume = {56},
  pages = {100900},
  issn = {13675788},
  doi = {10.1016/j.arcontrol.2023.100900},
  urldate = {2026-01-08},
  langid = {english}
}

@article{LopezNaranjoSalazarEtAl2025JComputPhys,
  title = {Weak-Form Modified Sparse Identification of Nonlinear Dynamics},
  author = {L{\'o}pez, Cristian and Naranjo, {\'A}ngel and Salazar, Diego and Moore, Keegan J.},
  year = 2025,
  month = sep,
  journal = {Journal of Computational Physics},
  pages = {114410},
  issn = {00219991},
  doi = {10.1016/j.jcp.2025.114410},
  urldate = {2025-10-04},
  langid = {english}
}

@article{LoreDePascualeLaiuEtAl2023NuclFusion,
  title = {Time-Dependent {{SOLPS-ITER}} Simulations of the Tokamak Plasma Boundary for Model Predictive Control Using {{SINDy}}},
  author = {Lore, J.D. and De Pascuale, S. and Laiu, P. and Russo, B. and Park, J.-S. and Park, J.M. and Brunton, S.L. and Kutz, J.N. and Kaptanoglu, A.A.},
  year = 2023,
  month = apr,
  journal = {Nuclear Fusion},
  volume = {63},
  number = {4},
  pages = {046015},
  issn = {0029-5515, 1741-4326},
  doi = {10.1088/1741-4326/acbe0e},
  urldate = {2023-03-27},
  langid = {english}
}

@article{Lorenz1963JAtmosSci,
  title = {Deterministic {{Nonperiodic Flow}}},
  author = {Lorenz, Edward N.},
  year = 1963,
  month = mar,
  journal = {Journal of the Atmospheric Sciences},
  volume = {20},
  number = {2},
  pages = {130--141},
  issn = {0022-4928, 1520-0469},
  doi = {10.1175/1520-0469(1963)020<0130:DNF>2.0.CO;2},
  urldate = {2024-08-24},
  langid = {english}
}

@article{ManganBruntonProctorEtAl2016IEEETransMolBiolMulti-ScaleCommun,
  title = {Inferring {{Biological Networks}} by {{Sparse Identification}} of {{Nonlinear Dynamics}}},
  author = {Mangan, Niall M. and Brunton, Steven L. and Proctor, Joshua L. and Kutz, J. Nathan},
  year = 2016,
  month = jun,
  journal = {IEEE Transactions on Molecular, Biological and Multi-Scale Communications},
  volume = {2},
  number = {1},
  pages = {52--63},
  issn = {2372-2061, 2332-7804},
  doi = {10.1109/TMBMC.2016.2633265},
  urldate = {2019-10-27},
  langid = {english}
}

@article{ManzoorPeiSunEtAl2022Drones,
  title = {Model {{Predictive Control Technique}} for {{Ducted Fan Aerial Vehicles Using Physics-Informed Machine Learning}}},
  author = {Manzoor, Tayyab and Pei, Hailong and Sun, Zhongqi and Cheng, Zihuan},
  year = 2022,
  month = dec,
  journal = {Drones},
  volume = {7},
  number = {1},
  pages = {4},
  issn = {2504-446X},
  doi = {10.3390/drones7010004},
  urldate = {2026-04-07},
  langid = {english}
}

@article{MayneRawlingsRaoEtAl2000Automatica,
  title = {Constrained Model Predictive Control: {{Stability}} and Optimality},
  author = {Mayne, D Q and Rawlings, J B and Rao, C V and Scokaert, P O M},
  year = 2000,
  journal = {Automatica},
  volume = {36},
  number = {6},
  pages = {789--814},
  doi = {10.1016/S0005-1098(99)00214-9},
  langid = {english}
}

@article{MessengerBortz2021JComputPhys,
  title = {Weak {{SINDy For Partial Differential Equations}}},
  author = {Messenger, Daniel A. and Bortz, David M.},
  year = 2021,
  month = oct,
  journal = {Journal of Computational Physics},
  volume = {443},
  pages = {110525},
  doi = {10.1016/j.jcp.2021.110525},
  urldate = {2020-07-07}
}

@article{MessengerBortz2021MultiscaleModelSimul,
  title = {Weak {{SINDy}}: {{Galerkin-Based Data-Driven Model Selection}}},
  shorttitle = {Weak {{SINDy}}},
  author = {Messenger, Daniel A. and Bortz, David M.},
  year = 2021,
  journal = {Multiscale Modeling \& Simulation},
  volume = {19},
  number = {3},
  pages = {1474--1497},
  doi = {10.1137/20M1343166},
  urldate = {2020-07-07}
}

@article{MessengerBortz2022PhysicaD,
  title = {Learning Mean-Field Equations from Particle Data Using {{WSINDy}}},
  author = {Messenger, Daniel A. and Bortz, David M.},
  year = 2022,
  month = nov,
  journal = {Physica D: Nonlinear Phenomena},
  volume = {439},
  pages = {133406},
  issn = {01672789},
  doi = {10.1016/j.physd.2022.133406},
  urldate = {2022-07-20},
  langid = {english}
}

@article{MessengerDwyerDukic2024JRSocInterface,
  title = {Weak-Form Inference for Hybrid Dynamical Systems in Ecology},
  author = {Messenger, Daniel and Dwyer, Greg and Dukic, Vanja},
  year = 2024,
  month = dec,
  journal = {Journal of The Royal Society Interface},
  volume = {21},
  number = {221},
  pages = {20240376},
  issn = {1742-5662},
  doi = {10.1098/rsif.2024.0376},
  urldate = {2026-02-11},
  langid = {english}
}

@article{MichaelMellingerLindseyEtAl2010IEEERobotAutomatMag,
  title = {The {{GRASP Multiple Micro-UAV Testbed}}},
  author = {Michael, Nathan and Mellinger, Daniel and Lindsey, Quentin and Kumar, Vijay},
  year = 2010,
  month = sep,
  journal = {IEEE Robotics \& Automation Magazine},
  volume = {17},
  number = {3},
  pages = {56--65},
  issn = {1070-9932},
  doi = {10.1109/MRA.2010.937855},
  urldate = {2026-01-08},
  copyright = {https://ieeexplore.ieee.org/Xplorehelp/downloads/license-information/IEEE.html},
  langid = {english}
}

@article{MinorMessengerDukicEtAl2025JournalofGeophysicalResearchMachineLearningandComputation,
  title = {Learning {{Physically Interpretable Atmospheric Models From Data With WSINDy}}},
  author = {Minor, Seth and Messenger, Daniel A. and Dukic, Vanja and Bortz, David M.},
  year = 2025,
  month = sep,
  journal = {Journal of Geophysical Research: Machine Learning and Computation},
  volume = {2},
  number = {3},
  pages = {e2025JH000602},
  issn = {2993-5210, 2993-5210},
  doi = {10.1029/2025JH000602},
  urldate = {2026-02-11},
  langid = {english}
}

@article{OsmanXiaMahdiEtAl2025IntlJRobustNonlinear,
  title = {An {{Adaptive SINDy}}-{{Lyapunov Model Predictive Control Framework}} for {{Dual}}-{{System VTOL UAVs}}},
  author = {Osman, Mohammed and Xia, Yuanqing and Mahdi, Mohammed and Manzoor, Tayyab and Bajodah, Abdulrahman H. and Ali, Asif and Ali, Abid and Ahmed, Azzam},
  year = 2025,
  month = nov,
  journal = {International Journal of Robust and Nonlinear Control},
  pages = {rnc.70272},
  issn = {1049-8923, 1099-1239},
  doi = {10.1002/rnc.70272},
  urldate = {2026-01-08},
  langid = {english}
}

@article{ProctorBruntonKutz2016SIAMJApplDynSyst,
  title = {Dynamic {{Mode Decomposition}} with {{Control}}},
  author = {Proctor, Joshua L. and Brunton, Steven L. and Kutz, J. Nathan},
  year = 2016,
  month = jan,
  journal = {SIAM Journal on Applied Dynamical Systems},
  volume = {15},
  number = {1},
  pages = {142--161},
  issn = {1536-0040},
  doi = {10.1137/15M1013857},
  urldate = {2019-06-30},
  langid = {english}
}

@article{RaissiKarniadakis2018JournalofComputationalPhysics,
  title = {Hidden Physics Models: {{Machine}} Learning of Nonlinear Partial Differential Equations},
  shorttitle = {Hidden Physics Models},
  author = {Raissi, Maziar and Karniadakis, George Em},
  year = 2018,
  month = mar,
  journal = {Journal of Computational Physics},
  volume = {357},
  pages = {125--141},
  issn = {00219991},
  doi = {10.1016/j.jcp.2017.11.039},
  urldate = {2026-02-21},
  langid = {english}
}

@book{RawlingsMayneDiehl2024,
  title = {Model {{Predictive Control}}: {{Theory}}, {{Computation}}, and {{Design}}},
  author = {Rawlings, J.B. and Mayne, D.Q. and Diehl, M.M.},
  year = 2024,
  edition = {2nd},
  publisher = {Nob Hill Publishing, LLC},
  address = {Santa Barbara, CA}
}

@article{ReinboldGrigoriev2019PhysRevE,
  title = {Data-Driven Discovery of Partial Differential Equation Models with Latent Variables},
  author = {Reinbold, Patrick A. K. and Grigoriev, Roman O.},
  year = 2019,
  month = aug,
  journal = {Physical Review E},
  volume = {100},
  number = {2},
  pages = {022219},
  issn = {2470-0045, 2470-0053},
  doi = {10.1103/PhysRevE.100.022219},
  urldate = {2021-11-10},
  langid = {english}
}

@article{ReiterBaelmansBorner2005FusionScienceandTechnology,
  title = {The {{EIRENE}} and {{B2-EIRENE Codes}}},
  author = {Reiter, D. and Baelmans, M. and B{\"o}rner, P.},
  year = 2005,
  month = feb,
  journal = {Fusion Science and Technology},
  volume = {47},
  number = {2},
  pages = {172--186},
  issn = {1536-1055, 1943-7641},
  doi = {10.13182/FST47-172},
  urldate = {2026-01-08},
  langid = {english}
}

@article{RummelMessengerBeckerEtAl2025arXiv250208881,
  title = {{{WENDy}} for {{Nonlinear-in-Parameter ODEs}}},
  author = {Rummel, Nicholas and Messenger, Daniel A. and Becker, Stephen and Dukic, Vanja and Bortz, David M.},
  year = 2025,
  month = feb,
  journal = {arXiv:2502.08881},
  eprint = {2502.08881},
  archiveprefix = {arXiv}
}

@article{SchaefferMcCalla2017PhysRevE,
  title = {Sparse Model Selection via Integral Terms},
  author = {Schaeffer, Hayden and McCalla, Scott G.},
  year = 2017,
  month = aug,
  journal = {Physical Review E},
  volume = {96},
  number = {2},
  pages = {023302},
  issn = {2470-0045, 2470-0053},
  doi = {10.1103/PhysRevE.96.023302},
  urldate = {2021-09-27},
  langid = {english}
}

@article{Schmid2010JFluidMech,
  title = {Dynamic Mode Decomposition of Numerical and Experimental Data},
  author = {Schmid, Peter J.},
  year = 2010,
  month = aug,
  journal = {Journal of Fluid Mechanics},
  volume = {656},
  pages = {5--28},
  issn = {0022-1120, 1469-7645},
  doi = {10.1017/S0022112010001217},
  urldate = {2019-07-02},
  langid = {english}
}

@article{SchmidDoostanPourahmadian2024arXiv240920510,
  title = {Ensemble {{WSINDy}} for {{Data Driven Discovery}} of {{Governing Equations}} from {{Laser-based Full-field Measurements}}},
  author = {Schmid, Abigail C. and Doostan, Alireza and Pourahmadian, Fatemeh},
  year = 2024,
  month = sep,
  journal = {arXiv:2409.20510},
  eprint = {2409.20510},
  primaryclass = {math},
  urldate = {2024-10-10},
  archiveprefix = {arXiv},
  langid = {english}
}

@article{SchmidtLipson2009Science,
  title = {Distilling {{Free-Form Natural Laws}} from {{Experimental Data}}},
  author = {Schmidt, Michael and Lipson, Hod},
  year = 2009,
  month = apr,
  journal = {Science},
  volume = {324},
  number = {5923},
  pages = {81--85},
  issn = {0036-8075, 1095-9203},
  doi = {10.1126/science.1165893},
  urldate = {2019-06-28},
  copyright = {American Association for the Advancement of Science},
  langid = {english},
  pmid = {19342586}
}

@article{Schroeder1970IEEETransInformTheory,
  title = {Synthesis of Low-Peak-Factor Signals and Binary Sequences with Low Autocorrelation ({{Corresp}}.)},
  author = {Schroeder, M.},
  year = 1970,
  month = jan,
  journal = {IEEE Transactions on Information Theory},
  volume = {16},
  number = {1},
  pages = {85--89},
  issn = {0018-9448, 1557-9654},
  doi = {10.1109/TIT.1970.1054411},
  urldate = {2025-09-10},
  copyright = {https://ieeexplore.ieee.org/Xplorehelp/downloads/license-information/IEEE.html},
  langid = {english}
}

@article{SchwenzerAyBergsEtAl2021IntJAdvManufTechnol,
  title = {Review on Model Predictive Control: An Engineering Perspective},
  shorttitle = {Review on Model Predictive Control},
  author = {Schwenzer, Max and Ay, Muzaffer and Bergs, Thomas and Abel, Dirk},
  year = 2021,
  month = nov,
  journal = {The International Journal of Advanced Manufacturing Technology},
  volume = {117},
  number = {5-6},
  pages = {1327--1349},
  issn = {0268-3768, 1433-3015},
  doi = {10.1007/s00170-021-07682-3},
  urldate = {2025-09-08},
  langid = {english}
}

@techreport{Shinbrot1954NACATN3288,
  title = {On the Analysis of Linear and Nonlinear Dynamical Systems for Transient-Response Data},
  author = {Shinbrot, Marvin},
  year = 1954,
  month = dec,
  number = {NACA TN 3288},
  address = {Moffett Field, CA},
  institution = {Ames Aeronautical Laboratory}
}

@article{ShiriParkBennis2020IEEETransCommun,
  title = {Communication-{{Efficient Massive UAV Online Path Control}}: {{Federated Learning Meets Mean-Field Game Theory}}},
  shorttitle = {Communication-{{Efficient Massive UAV Online Path Control}}},
  author = {Shiri, Hamid and Park, Jihong and Bennis, Mehdi},
  year = 2020,
  month = nov,
  journal = {IEEE Transactions on Communications},
  volume = {68},
  number = {11},
  pages = {6840--6857},
  issn = {0090-6778, 1558-0857},
  doi = {10.1109/TCOMM.2020.3017281},
  urldate = {2026-06-10},
  copyright = {https://ieeexplore.ieee.org/Xplorehelp/downloads/license-information/IEEE.html},
  langid = {english}
}

@article{SorokinaSygletosTuritsyn2016OptExpress,
  title = {Sparse Identification for Nonlinear Optical Communication Systems: {{SINO}} Method},
  shorttitle = {Sparse Identification for Nonlinear Optical Communication Systems},
  author = {Sorokina, Mariia and Sygletos, Stylianos and Turitsyn, Sergei},
  year = 2016,
  month = dec,
  journal = {Optics Express},
  volume = {24},
  number = {26},
  pages = {30433},
  issn = {1094-4087},
  doi = {10.1364/OE.24.030433},
  urldate = {2025-09-30},
  copyright = {https://doi.org/10.1364/OA\_License\_v1\#VOR-OA},
  langid = {english}
}

@book{Strogatz2015,
  title = {Nonlinear Dynamics and Chaos: With Applications to Physics, Biology, Chemistry, and Engineering},
  shorttitle = {Nonlinear Dynamics and Chaos},
  author = {Strogatz, Steven H.},
  year = 2015,
  edition = {2nd},
  publisher = {CRC Press.},
  doi = {10.1201/9780429492563},
  isbn = {978-0-8133-4910-7},
  langid = {english},
  lccn = {Q172.5.C45 S767 2015}
}

@article{TranBortz2025arXiv250703206,
  title = {Weak {{Form Scientific Machine Learning}}: {{Test Function Construction}} for {{System Identification}}},
  author = {Tran, April and Bortz, David M.},
  year = 2025,
  month = jul,
  journal = {arXiv:2507.03206},
  eprint = {2507.03206},
  archiveprefix = {arXiv}
}

@article{LinHorneTinoEtAl1996IEEETransNeuralNetw,
  title = {Learning Long-Term Dependencies in {{NARX}} Recurrent Neural Networks},
  author = {Lin, Tsungnan and Horne, Bill G. and Tino, Peter and Giles, C. Lee},
  year = 1996,
  month = nov,
  journal = {IEEE Transactions on Neural Networks},
  volume = {7},
  number = {6},
  pages = {1329--1338},
  issn = {1045-9227, 1941-0093},
  doi = {10.1109/72.548162},
  urldate = {2026-04-09},
  copyright = {https://ieeexplore.ieee.org/Xplorehelp/downloads/license-information/IEEE.html},
  langid = {english}
}

@article{TuRowleyLuchtenburgEtAl2014JComputDyn,
  title = {On Dynamic Mode Decomposition: {{Theory}} and Applications},
  shorttitle = {On Dynamic Mode Decomposition},
  author = {Tu, Jonathan H. and Rowley, Clarence W. and Luchtenburg, Dirk M. and Brunton, Steven L. and Kutz, J. Nathan},
  year = 2014,
  journal = {Journal of Computational Dynamics},
  volume = {1},
  number = {2},
  pages = {391--421},
  issn = {2158-2505},
  doi = {10.3934/jcd.2014.1.391},
  urldate = {2025-09-30},
  langid = {english}
}

@article{VaseyMessengerBortzEtAl2025JComputPhys,
  title = {Influence of Initial Conditions on Data-Driven Model Identification and Information Entropy for Ideal Mhd Problems},
  author = {Vasey, Gina and Messenger, Daniel A. and Bortz, David M. and Christlieb, Andrew and O'Shea, Brian},
  year = 2025,
  month = mar,
  journal = {Journal of Computational Physics},
  volume = {524},
  pages = {113719},
  doi = {10.1016/j.jcp.2025.113719},
  urldate = {2023-12-14},
  langid = {english}
}

@article{WangLiWenEtAl2026IEEETransVehTechnol,
  title = {A {{Structured Pruning}} Based {{MPC}} for {{Virtual Coupling Trains}} via {{Utilizing WHT-SINDy}}},
  author = {Wang, Jincheng and Li, Kun and Wen, Tao and Roberts, Clive},
  year = 2026,
  journal = {IEEE Transactions on Vehicular Technology},
  pages = {1--15},
  issn = {0018-9545, 1939-9359},
  doi = {10.1109/TVT.2026.3701305},
  urldate = {2026-06-14},
  copyright = {https://ieeexplore.ieee.org/Xplorehelp/downloads/license-information/IEEE.html},
  langid = {english}
}

@article{WangYangLaiEtAl2011PhysRevLett,
  title = {Predicting {{Catastrophes}} in {{Nonlinear Dynamical Systems}} by {{Compressive Sensing}}},
  author = {Wang, Wen-Xu and Yang, Rui and Lai, Ying-Cheng and Kovanis, Vassilios and Grebogi, Celso},
  year = 2011,
  month = apr,
  journal = {Physical Review Letters},
  volume = {106},
  number = {15},
  pages = {154101},
  issn = {0031-9007, 1079-7114},
  doi = {10.1103/PhysRevLett.106.154101},
  urldate = {2020-09-07},
  langid = {english}
}

@article{WentzDoostan2023ComputerMethodsinAppliedMechanicsandEngineering,
  title = {Derivative-Based {{SINDy}} ({{DSINDy}}): {{Addressing}} the Challenge of Discovering Governing Equations from Noisy Data},
  shorttitle = {Derivative-Based {{SINDy}} ({{DSINDy}})},
  author = {Wentz, Jacqueline and Doostan, Alireza},
  year = 2023,
  month = aug,
  journal = {Computer Methods in Applied Mechanics and Engineering},
  volume = {413},
  pages = {116096},
  issn = {00457825},
  doi = {10.1016/j.cma.2023.116096},
  urldate = {2023-11-19},
  langid = {english}
}

@article{WiltingPriesemann2018NatCommun,
  title = {Inferring Collective Dynamical States from Widely Unobserved Systems},
  author = {Wilting, Jens and Priesemann, Viola},
  year = 2018,
  month = jun,
  journal = {Nature Communications},
  volume = {9},
  number = {1},
  pages = {2325},
  issn = {2041-1723},
  doi = {10.1038/s41467-018-04725-4},
  urldate = {2025-09-10},
  langid = {english}
}

@article{WoodallEsparzaGutovaEtAl2024APLBioeng,
  title = {Model Discovery Approach Enables Noninvasive Measurement of Intra-Tumoral Fluid Transport in Dynamic {{MRI}}},
  author = {Woodall, Ryan T. and Esparza, Cora C. and Gutova, Margarita and Wang, Maosen and Cunningham, Jessica J. and Brummer, Alexander B. and Stine, Caleb A. and Brown, Christine C. and Munson, Jennifer M. and Rockne, Russell C.},
  year = 2024,
  month = jun,
  journal = {APL Bioengineering},
  volume = {8},
  number = {2},
  pages = {026106},
  issn = {2473-2877},
  doi = {10.1063/5.0190561},
  urldate = {2024-06-12},
  langid = {english}
}

@article{WuTranRinconEtAl2019AIChEJournal,
  title = {Machine Learning-based Predictive Control of Nonlinear Processes. {{Part I}}: {{Theory}}},
  shorttitle = {Machine Learning-based Predictive Control of Nonlinear Processes. {{Part I}}},
  author = {Wu, Zhe and Tran, Anh and Rincon, David and Christofides, Panagiotis D.},
  year = 2019,
  month = nov,
  journal = {AIChE Journal},
  volume = {65},
  number = {11},
  pages = {e16729},
  issn = {0001-1541, 1547-5905},
  doi = {10.1002/aic.16729},
  urldate = {2026-06-10},
  langid = {english}
}

@article{WuTranRinconEtAl2019AIChEJournala,
  title = {Machine-learning-based Predictive Control of Nonlinear Processes. {{Part II}}: {{Computational}} Implementation},
  shorttitle = {Machine-learning-based Predictive Control of Nonlinear Processes. {{Part II}}},
  author = {Wu, Zhe and Tran, Anh and Rincon, David and Christofides, Panagiotis D.},
  year = 2019,
  month = nov,
  journal = {AIChE Journal},
  volume = {65},
  number = {11},
  pages = {e16734},
  issn = {0001-1541, 1547-5905},
  doi = {10.1002/aic.16734},
  urldate = {2026-06-01},
  langid = {english}
}

@article{WuChristofidesWuEtAl2025RevChemEng,
  title = {A Tutorial Review of Machine Learning-Based Model Predictive Control Methods},
  author = {Wu, Zhe and Christofides, Panagiotis D. and Wu, Wanlu and Wang, Yujia and Abdullah, Fahim and Alnajdi, Aisha and Kadakia, Yash},
  year = 2025,
  month = may,
  journal = {Reviews in Chemical Engineering},
  volume = {41},
  number = {4},
  pages = {359--400},
  issn = {0167-8299, 2191-0235},
  doi = {10.1515/revce-2024-0055},
  urldate = {2026-06-01},
  copyright = {http://creativecommons.org/licenses/by/4.0},
  langid = {english}
}

@article{YahagiSetoYonezawaEtAl2025IntJControlAutomSyst,
  title = {Sparse {{Identification}} and {{Nonlinear Model Predictive Control}} for {{Diesel Engine Air Path System}}},
  author = {Yahagi, Shuichi and Seto, Hiroki and Yonezawa, Ansei and Kajiwara, Itsuro},
  year = 2025,
  month = feb,
  journal = {International Journal of Control, Automation and Systems},
  volume = {23},
  number = {2},
  pages = {620--629},
  issn = {2005-4092},
  doi = {10.1007/s12555-024-0452-9},
  urldate = {2025-02-18},
  langid = {english}
}

@article{YanWang2012IEEETransIndInf,
  title = {Model {{Predictive Control}} of {{Nonlinear Systems With Unmodeled Dynamics Based}} on {{Feedforward}} and {{Recurrent Neural Networks}}},
  author = {Yan, Zheng and Wang, Jun},
  year = 2012,
  month = nov,
  journal = {IEEE Transactions on Industrial Informatics},
  volume = {8},
  number = {4},
  pages = {746--756},
  issn = {1551-3203, 1941-0050},
  doi = {10.1109/TII.2012.2205582},
  urldate = {2025-10-11},
  langid = {english}
}

@misc{YazdannikTayefiFarrokh2023arXiv231112906,
  title = {Nonlinear {{System Identification}} of {{Swarm}} of {{UAVs Using Deep Learning Methods}}},
  author = {Yazdannik, Saman and Tayefi, Morteza and Farrokh, Mojtaba},
  year = 2023,
  month = nov,
  number = {arXiv:2311.12906},
  eprint = {2311.12906},
  publisher = {arXiv},
  doi = {10.48550/arXiv.2311.12906},
  urldate = {2026-06-10},
  archiveprefix = {arXiv},
  langid = {english}
}

@misc{YirenLu2022,
  title = {Yrlu/Quadrotor: {{Quadrotor Control}}, {{Path Planning}} and {{Trajectory Optimization}}},
  shorttitle = {Yrlu/Quadrotor},
  author = {{Yiren Lu}},
  year = 2022,
  month = jul,
  doi = {10.5281/ZENODO.6796214},
  urldate = {2026-02-20},
  copyright = {Open Access},
  howpublished = {Zenodo}
}

@article{ZhangGuoZhangEtAl2024PhilTransRSocA,
  title = {Physics-Informed Deep Learning for Structural Vibration Identification and Its Application on a Benchmark Structure},
  author = {Zhang, Minte and Guo, Tong and Zhang, Guodong and Liu, Zhongxiang and Xu, Weijie},
  year = 2024,
  month = jan,
  journal = {Philosophical Transactions of the Royal Society A: Mathematical, Physical and Engineering Sciences},
  volume = {382},
  number = {2264},
  pages = {20220400},
  issn = {1364-503X, 1471-2962},
  doi = {10.1098/rsta.2022.0400},
  urldate = {2026-01-09},
  langid = {english}
}
\bibliographystyle{RS}
\end{document}